\documentclass[10pt,reqno,oneside]{amsproc}
\usepackage[cp1251]{inputenc}
\usepackage[russian]{babel}
\usepackage{amssymb,epic,eepic}
\allowdisplaybreaks

\title[]{Поведение асимптотики положительного спектра семейства периодических
задач Штурма--Лиувилля при непрерывном переходе от дефинитной к индефинитной
задаче}
\author{Д.~А.~Попов}

\renewcommand{\Re}{\operatorname{Re}}
\renewcommand{\Im}{\operatorname{Im}}
\makeatletter
\@addtoreset{equation}{section}
\makeatother

\newtheorem*{theor*}{Теорема}
\newtheorem*{los}{Лемма о спектре}
\begin{document}
\begin{abstract}
Рассматривается задача о спектре зависящего от параметра \(a\in\mathbb R\)
семейства периодических задач Штурма--Лиувилля для уравнения вида
\(u''+\lambda^2(g(x)-a)\,u=0\), где \(\lambda\) "--- спектральный параметр.
Предполагается, что \(g:\mathbb R\to\mathbb R\) "--- достаточно гладкая
периодическая функция, имеющая на периоде один простой максимум \(g(x_{max})=
a_1>0\) и один простой минимум \(g(x_{min})=a_2>0\). Кроме того,
предполагается, что функции \(g(x-x_{min})\) и \(g(x-x_{min})\) являются
чётными. При этих предположениях на всём интервале \(0\leqslant a<a_1\),
включая окрестности точек \(a=a_1\) и \(a=a_2\), явно вычислены первые два
члена асимптотики положительных собственных значений. Показано, что при
\(\lambda\gg 1\) спектр состоит из двух ветвей \(\lambda=\lambda_{\pm}(a,p)\),
нумеруемых выбором знака \(\pm\) и целым числом \(p\in\mathbb Z^+\),
\(p\gg 1\). Получена единая интерполяционная формула, описывающая поведение
асимптотики ветвей спектра при переходе от дефинитной (классической) при
\(a<a_2\) к индефинитной при \(a>a_2\) задаче.
\end{abstract}
\begin{flushleft}
УДК~517.927.25
\end{flushleft}
\maketitle
\section{Введение}\label{pt:1}
В настоящей работе рассматривается периодическая задача Штурма--Лиувилля
\begin{equation}\label{eq:1:1}
	\begin{gathered}
		\dfrac{d^2u}{dx^2}(x)+\lambda^2(g(x)-a)\,u(x)=0,\\
		u(-\pi)=u(\pi),\qquad \dfrac{du}{dx}(-\pi)=\dfrac{du}{dx}(\pi).
	\end{gathered}
\end{equation}
Здесь \(g:\mathbb R\to\mathbb R\) "--- достаточно гладкая периодическая,
с периодом \(2\pi\), функция, \(\lambda\) "--- спектральный параметр, и
\(a\in\mathbb R\) "--- параметр. Функцию \(g\) можно рассматривать как функцию
на окружности \(\mathbb S^1\), причём предполагается, что там она имеет
один простой максимум \(g(x_{max})=a_1\) и один простой минимум \(g(x_{min})=
a_2\). Точное описание класса \(\mathcal G\) рассматриваемых функций \(g\)
будет дано в~\ref{pt:2}. Уравнение~\eqref{eq:1:1} не меняется при
одновременной замене \(a\mapsto a+C\), \(g(x)\mapsto g(x)+C\), и поэтому мы
будем считать, что \(0\leqslant a_2\leqslant g(x)\leqslant a_1\). Так как
при \(a>a_1\) уравнение~\eqref{eq:1:1} не имеет периодических решений, то
ограничимся рассмотрением интервала \(U\), в котором \(0<a<a_1\).
В рассматриваемой задаче при \(0\leqslant a<a_2\) нет точек поворота "---
нулей функции \(r(x)=g(x)-a\). При \(a_2<a<a_1\) таких точек две, и они
сливаются при \(a=a_i\) (\(i=1,2\)). Таким образом, при \(a=a_2\) происходит
непрерывный переход от дефинитной (классической) к индефинитной при \(a>a_2\)
задаче Штурма--Лиувилля. Напомним, что краевая задача для уравнения
\(-y''+q(x)y=\mu\,r(x)y\) (\(a\leqslant x\leqslant b\)) называется дефинитной,
если \(r(x)>0\) (\(x\in [a,b]\)), и индефинитной, если функция \(r(x)\) меняет
знак на интервале \((a,b)\).

Цель настоящей работы состоит в построении равномерной по параметру \(a\in U\)
на всём полуинтервале \(U\) (\(0\leqslant a<a_1\)) асимптотики положительных
собственных значений \(\lambda=\lambda(a)\) при \(\lambda\to\infty\). Ниже такая
асимптотика будет построена при предположении чётности функции \(g(x)\),
в соответствии с которым функции \(g(x-x_{min})\) и \(g(x-x_{max})\)
предполагаются чётными.

При \(a<a_1\) функция \(r(x)=g(x)-a>0\), и с помощью преобразования
Лиувилля \cite{1,2} \(y(z)=(r(x))^{1/4}u(x)\), \(z=B^{-1}\int_{-x}^x
r^{1/2}(x)\,dx\), \(B=\int_{-\pi}^{\pi} r^{1/2}(x)\,dx\) задача~\eqref{eq:1:1}
сводится к классической задаче Штурма--Лиувилля
\begin{gather}\label{eq:1:2}
	\begin{gathered}
		\dfrac{d^2y}{dz^2}(z)+(\tilde\lambda-q(z))y(z)=0,\qquad
		(0<z<1),\\
		y(0)=y(1),\qquad \dfrac{dy}{dz}(0)=\dfrac{dy}{dz}(1).
	\end{gathered}
	\intertext{При этом}
	q(z)=\dfrac{B^2}{4r^2}\left(\dfrac{d^2r}{dx^2}-\dfrac{5}{4}\cdot
	\dfrac{1}{r}\cdot\dfrac{dr}{dx}\right),\qquad\tilde\lambda=B^2\lambda^2.
\end{gather}
Спектр задачи~\eqref{eq:1:2} хорошо известен \cite{3,4,5}. Собственные значения
\(\tilde\lambda_k\) образуют бесконечную последовательность вида
\begin{equation}\label{eq:1:4}
	-\infty<\tilde\lambda_0<\tilde\lambda_1\leqslant\tilde\lambda_2<\ldots
	<\tilde\lambda_{2n+3}\leqslant\tilde\lambda_{2n-2}<\tilde\lambda_{2n-1}
	\leqslant\tilde\lambda_{2n}<\tilde\lambda_{2n+1}\leqslant\ldots
\end{equation}
При этом собственные функции \(y_{2n-1}(z)\), \(y_{2n}(z)\), отвечающие
собственным значениям \(\tilde\lambda_{2n-1}\), \(\tilde\lambda_{2n}\),
имеют \(2n\) нулей на полуинтервале \([0,1)\), и \(\tilde\lambda_{2n-1}\),
\(\tilde\lambda_{2n}\) имеют одинаковое асимптотическое разложение вида
\begin{equation}\label{eq:1:5}
	\sqrt{\tilde\lambda_{2n}}=2\pi n+\dfrac{c_1}{n}+\dfrac{c_3}{n^3}+\ldots,
\end{equation}
если функция \(q(z)\) "--- гладкая. В интересующем нас случае \(q(z)\to\infty\)
при \(a\to a_2-0\), и разложение~\eqref{eq:1:5} не является равномерным по \(a\)
во всей области \(a\leqslant a_2\).

Рассмотрим область \(a>a_2\). В этой области функция \(r(x)=g(x)-a\) меняет
знак, и мы имеем дело с индефинитной задачей Штурма--Лиувилля. Эта задача
для уравнения \(y''+(\mu\,r-q)y=0\) на конечном интервале \((d_1,d_2)\)
при разделённых граничных условиях \(y(d_i)\,\cos\alpha_i-y'(d_i)\,
\sin\alpha_i=0\) (\(i=1,2\)) рассматривалась в ряде сравнительно недавних
работ \cite{6}--\cite{13}. С точки зрения спектральной теории, индефинитная
задача вписывается в рамки теории операторных пучков и теории пространств
с индефинитной метрикой. В указанных работах, в частности, показано, что
рассматриваемая задача может иметь конечное число невещественных собственных
значений, и существуют две бесконечные последовательности \(\mu^{\pm}_n\to
\pm\infty\) (\(n\to\infty\)) вещественных собственных значений, для которых
имеет место асимптотическая формула
\begin{equation}\label{eq:1:6}
	\mu^{\pm}_n\sim\pm\pi^2 n^2\left[\int\limits_{d_1}^{d_2}
	(r_{\pm}(x))^{1/2}\,dx\right]^{-2},\qquad
	r_{\pm}(x)=\max\{\pm r(x),0\}.
\end{equation}
Полное асимптотическое разложение вещественных собственных значений для задачи
с разделёнными граничными условиями и \(N\) простыми точками поворота получено
в работе~\cite{14}. Периодическая задача при аналогичных предположениях
рассматривалась в работе~\cite{15}. Метод склеивания асимптотических
разложений, используемый в этих работах, предполагает, что расстояние между
точками поворота ограничено снизу, и не позволяет получить равномерные
асимптотики, пригодные в окрестностях точек \(a=a_i\) (\(i=1,2\)).

Так как при выполнении условия чётности периодическая задача сводится к
двум задачам с разделёнными граничными условиями, то в области \(a_1-c_1
\leqslant a\leqslant a_2+c_2\) (\(c_i>0\)) в нашей задаче можно
воспользоваться результатами работы~\cite{14}. В настоящей работе метод
склеивания не используется и используемый ниже метод модельных уравнений
позволяет исследовать асимптотику спектра при всех \(a\in U\).

Уже давно было ясно, что при возможности слияния двух простых точек поворота
в качестве модельного уравнения необходимо использовать уравнение Вебера
(см.~книги \cite{1,2}, работу \cite{16} и цитируемую там литературу).
Исследованию спектра задачи~\eqref{eq:1:1} на основе этой идеи посвящена
первая часть работы \cite{16}. Результат работы \cite{16}, касающийся
рассматриваемой задачи, можно сформулировать следующим образом. При \(\lambda
\gg 1\) имеются две ветви спектра \(\lambda_{\pm}(a,p)\), где \(p\in
\mathbb Z^+\), \(p\gg 1\), и они находятся из уравнений
\begin{gather}\label{eq:1:7}
	\lambda F(a)=2\pi p+G_{\pm}(\lambda,a).\\
	\intertext{В этих уравнениях}\label{eq:1:8}
	F(a)=\int\limits_{g(x)>a}\sqrt{g(x)-a}\,dx.
\end{gather}
В формуле~\eqref{eq:1:8} интегрирование ведётся по пересечению области
\(g(x)>a\) с периодом и
\begin{equation}\label{eq:1:9}
	G_{\pm}(\lambda,a)=\left\{\begin{array}{l@{\qquad}l}
	\pm\dfrac{\pi}{2}+O(\lambda^{-2/3}\ln\lambda)&
	(a_2+C\lambda^{-2/3}\leqslant a<a_1)\\
	O(\lambda^{-2/3}\ln\lambda)& (0\leqslant a\leqslant a_2-C\lambda^{-2/3}),
	\end{array}\right.
\end{equation}
где \(C>0\) "--- некоторая константа. Используемый метод позволяет найти
асимптотику функций \(G_{\pm}(\lambda,a)\) и при \(|a-a_2|\leqslant
C\lambda^{-2/3}\). В работе \cite{16} утверждается, что асимптотически весь
положительный спектр сводится к двум указанным ветвям \(\lambda_{\pm}(a,p)\).
Это находится в полном соответствии с результатами работ \cite{6}--\cite{9},
\cite{14}, в которых рассматривается задача с разделёнными граничными условиями.

Заметим, что хотя точки спектра \(\lambda_{\pm}(a,p)\) нумеруются целым числом
\(p\), но остаётся задача о связи числа \(p\) с естественным номером \(n\)
собственных значений \(\lambda_{n,a}\), упорядоченных по величине
\[
	\ldots\leqslant\lambda_{n-1,a}\leqslant\lambda_{n,a}\leqslant
	\lambda_{n+1,a}\leqslant\ldots
\]
При решении этой задачи надо учитывать, что для индефинитной периодической
задачи осцилляционная теорема, устанавливающая связь номера \(n\) с числом
нулей соответствующей собственной функции, по-видимому, до сих пор
не доказана.

Так как стандартному ВКБ приближению отвечает \(G_{\pm}=0\) при \(a<a_2\)
и \(G_{\pm}=\pm\pi/2\) при \(a>a_2\), то вышеприведённый результат~\eqref{eq:1:7},
\eqref{eq:1:9} состоит в указании границ применимости этого приближения
и его точности. К сожалению, неточности, допущенные в работе \cite{16},
приводят к тому, что эти границы указаны неверно, и равенство~\eqref{eq:1:9}
неверно при \(|a-a_2|<C\lambda^{-1/3}\). Заметим, что результат~\eqref{eq:1:9}
противоречит результатам работ \cite{14,15}, из которых, так же как и из общей
структуры ВКБ разложения, следует, что в выражении \(G_{\pm}(\lambda,a)\) должен
содержаться член вида \(C(\lambda|a-a_2|)^{-1}\) (см.~\eqref{eq:1:10}). Этот
член при \(|a-a_2|\leqslant\lambda^{-2/3}\) имеет порядок \(\lambda^{-1/3}\),
что противоречит каждому из равенств~\eqref{eq:1:9}.

Для получения равномерных асимптотик спектра в настоящей работе используется
тот же метод, что и в работе \cite{16}. Этот метод состоит в использовании
полученных Олвером в работе \cite{17} равномерных по \(a\) асимптотик
фундаментальной системы решений уравнения~\eqref{eq:1:1} с равномерными оценками
остаточных членов в этих асимптотиках. Класс \(\mathcal G\) рассматриваемых
ниже функций по существу совпадает с классом функций, рассматриваемых Олвером.
Отметим при этом, что в работе Олвера условие "`чётности"' предполагается
выполненным только в области \(a<a_2\). Таким образом, ниже речь идёт
об аккуратном вычислении спектра в задаче~\eqref{eq:1:1} на основе результатов
Олвера. В частности будет показано, что в формуле~\eqref{eq:1:7}
(см.~\eqref{eq:2:13}, \eqref{eq:2:16})
\begin{equation}\label{eq:1:10}
	G_{\pm}(\lambda,a)=\left\{\begin{array}{l@{\qquad}l}
	\pm\dfrac{\pi}{2}-\dfrac{1}{12\lambda\alpha_2^2}+O(\lambda^{-2/3}
	\ln\lambda)&(a-a_2>C\lambda^{-7/9}(\ln\lambda)^{-1/3}),\\
	\dfrac{1}{12\lambda\alpha_2^2}+O(\lambda^{-1/2}(\ln\lambda)^{-1/2})&
	(a_2-a>C\lambda^{-2/3}(\ln\lambda)^{-1/3}).
	\end{array}\right.
\end{equation}
Величина \(\alpha_2\equiv\alpha_2(a)\) будет определена ниже
(см.~\eqref{eq:2:1}) и \(\alpha_2^2\asymp |a-a_2|\). Будет показано,
что во всём интервале \(U(0\leqslant a<a_1)\) имеет место
единая интерполяционная формула
\[
	G_{\pm}(\lambda,a)=H_{\pm}(b_2)+O(\lambda^{-1/2+\varepsilon})\qquad
	(\forall\varepsilon>0)
\]
с явным указанием вида функций \(H_{\pm}(x)\), и в этой формуле \(b_2=
(1/2)\cdot\lambda\alpha_2^2\) при \(a>a_2\) и \(b_2=-(1/2)\cdot\lambda
\alpha_2^2\) при \(a<a_2\). Заметим, что если \(g\) "--- гладкая функция,
то в области \(a_2-a>C>0\) выполняется соотношение \(|\lambda_+(a,p)-
\lambda_-(a,p)|=O(p^{\infty})\), а если \(g\) "--- целая функция, то,
как показано в работах \cite{18,19}, в этой области при \(p\to\infty\)
выполняется соотношение \(|\lambda_+(a,p)-\lambda_-(a,p)|\sim
e^{-\varkappa p}\).

Результаты Олвера позволяют (при \(g\in\mathcal G\)) получить равномерную
по параметру \(a\) асимптотику спектра и в задаче с разделёнными граничными
условиями. Краткие замечания по этому поводу сделаны в конце~\ref{pt:6}.
Выбор периодической задачи связан с указанным в работах \cite{15,16,20,21}
применением результатов в задаче о спектре оператора Лапласа--Бельтрами
на двумерном торе с метрикой Лиувилля.

Везде ниже предполагаются выполненными и больше не оговариваются следующие
соглашения о константах.

Буквой \(C\) обозначаются различные положительные константы, не зависящие
от параметра \(a\) и, возможно, зависящие от функции \(g\in\mathcal G\).
Все эти константы могут быть явно указаны.

Если утверждается, что для некоторой величины имеется оценка вида
\(O(f(\lambda,a))\) (\(O(f(p,a))\)), то это означает, что существуют
зависящие от только от функции \(g\in\mathcal G\) величины \(\Lambda\),
\(P\) и \(A\) такие, что \(\forall\lambda>\Lambda\) (\(\forall p>P\))
соответствующая величина по модулю меньше, чем \(Af(\lambda,a)\)
(\(Af(p,a)\)). Условия \(\lambda\gg 1\) (\(p\gg 1\)) означают, что
существуют величины \(\Lambda\) и \(P\), зависящие только от \(g\in
\mathcal G\) и такие, что соответствующие утверждения справедливы
\(\forall\lambda>\Lambda\) (\(\forall p>P\)). В большинстве случаев,
если не во всех, величины \(\Lambda\), \(P\) и \(A\) могут быть явно
указаны, но ниже этот вопрос не рассматривается.

Автор благодарен А.~А.~Шкаликову за стимулирующий интерес к этой работе.

\section{Формулировка результатов}\label{pt:2}
Введём класс \(\mathcal G\) рассматриваемых функций \(g\). Будем говорить,
что \(g\in\mathcal G\), если выполняются следующие условия (1)--(4).
\begin{enumerate}
\item \(g:\mathbb R\to\mathbb R\) "--- периодическая функция с периодом
\(2\pi\), имеющая шесть ограниченных производных \(g^{(n)}(x)\),
\(|g^{(n)}(x)|\leqslant C\) (\(n=0,1,\ldots,6\)).
\item Функция \(g\) имеет на полуинтервале \([-\pi,\pi)\) один простой
минимум в точке \(x_{min}\) и один простой максимум в точке \(x_{max}\).
При этом \(|g^{(1)}(x)|\neq 0\) (\(x\neq x_{min}, x_{max}\)),
\(g(x_{max})=a_1>0\), \(g(x_{min})=a_2>0\), \(g^{(2)}(x_{ex})\geqslant
C>0\) (\(x_{ex}=x_{min}, x_{max}\)).
\item Функции \(g(x-x_{min})\) и \(g(x-x_{max})\) "--- чётные.
\item Если \(x_{min}=0\), то при \(|x|\leqslant\pi\) выполняется тождество
\(g(x)=g_1(x^2)\), и чётная функция \(h(x)\), определённая при \(|x|\leqslant
\pi\) равенством \(h(x)=g_1(-x^2)\), имеет при \(|x|\leqslant\pi\) по крайней
мере шесть ограниченных производных \(|h^{(n)}(x)|\leqslant C\) (\(n=0,
1,\ldots,6\)). Имеется точка \(x_0\in [0,\pi-C]\), \(C>0\), такая, что
\(h(x_0)=0\). Функция \(h(x)\) при \(x\in[0,x_0]\) имеет один простой максимум
в точке \(x=0\), и \(h(0)=a_2\). При этом на отрезке \([0,x_0]\) функция
\(h(x)\) монотонно убывает, и
\[
	0\leqslant h(x)\leqslant a_2,\qquad h^{(1)}(x)<0,\qquad
	(0<|x|\leqslant x_0).
\]
\end{enumerate}

Условие~(3) будем называть условием чётности.

Заметим, что для аналитических функций \(g\) "--- \(h(x)=g(ix)\).

Приведём примеры функций из класса \(\mathcal G\). К классу \(\mathcal G\)
при \(x_{min}=0\) или \(x_{max}=0\) принадлежат периодические продолжения
функций вида \(g(x)=f(x^2)\) (\(0\leqslant x\leqslant\pi\)) при соответствующих
ограничениях на функцию \(f\). Можно показать, что при \(x_{min}=0\) функции
вида
\begin{gather*}
	g(x)=\sum\limits_{n=0}^{\infty} c_n\,\cos nx,\qquad
	|c_n|\leqslant Ce^{-\pi n(1+\delta)},\qquad (\delta>0)
	\intertext{при выполнении условий}
	c_n\leqslant 0\qquad (n\geqslant 1),\qquad
	|c_1|\leqslant\dfrac{\pi}{2}\sum\limits_{n=2}^{\infty}
	n\,|c_n|,\qquad\sum\limits_{n=1}^{\infty}|c_n|<c_0<
	\sum\limits_{n=1}^{\infty}|c_n|\,\ch(\pi n)
\end{gather*}
также принадлежат к классу \(\mathcal G\).

Интересующий нас полуинтервал \(0\leqslant a<a_1\) параметра \(a\)
обозначается через \(U\). Будем различать три области \(U_i\subset
U\) (\(i=1,2,3\)) параметра \(a\) такие, что \(U=\bigcup\limits_i
U_i\). По определению \(a\in U_1\), если \(a_2<a_0\leqslant a<a_1\);
\(a\in U_2\), если \(a_2\leqslant a\leqslant a'_0\) и \(a'_0\geqslant
a_0\); \(a\in U_3\), если \(0\leqslant a\leqslant a_2\). При \(a\in U_1\)
начало координат выбирается в точке \(x_{max}\) (\(x_{max}=0\)), а при
\(a\in U_2\cup U_3\) "--- в точке \(x_{min}\) (\(x_{min}=0\)). Таким
образом, ниже через \(g(x)\) в зависимости от \(a\) обозначаются
различные функции и, если при \(a\in U_1\) имеет место тождество
\(g(x)=\tilde g(x)\), то при \(a\in U_2\cup U_3\) имеет место тождество
\(g(x)=\tilde g(x-\pi)\).

Если \(a\in U_1\cup U_2\), то через \(x_1=-x_2\), \(x_2>0\) (\(x_i\equiv
x_i(a)\)) обозначаются точки поворота "--- корни уравнения \(g(x)=a\).
При \(a\in U_3\) через \(x_1=-x_2\), \(x_2>0\), обозначаются корни уравнения
\(h(x)=a\). Предполагается, что величины \(a_0\) и \(a'_0\) фиксированы
и выбраны таким образом, что \(\pi-x_2(a)\geqslant C>0\) при \(a=a_0,a'_0\).
Из условий~\eqref{eq:2:3} следует, что \(|x_{min}-x_{max}|=\pi\). Качественное
поведение функций \(g(x)\) и \(h(x)\) приведено на рис.~\ref{bil:1}.

\begin{figure}[ht]
\unitlength=0.132mm
\begin{picture}(650,600)
\put(0,100){\vector(1,0){500}}
\put(250,100){\vector(0,1){450}}
\put(505,110){{\(x\)}}
\put(260,555){{\(g(x)\)}}
\dashline[25]{25}(0,500)(500,500)
\dashline[25]{25}(140,350)(360,350)
\dashline[25]{25}(25,175)(475,175)
\dashline[25]{25}(25,100)(25,175)
\dashline[25]{25}(475,100)(475,175)
\dashline[25]{25}(150,100)(150,350)
\dashline[25]{25}(350,100)(350,350)
\put(25,70){{\(-\pi\)}}\put(475,70){{\(\pi\)}}
\put(150,70){{\(x_1\)}}\put(350,70){{\(x_2\)}}\put(250,70){{\(0\)}}
\put(250,100){\blacken\circle{10}}
\put(250,500){\blacken\circle{10}}
\put(250,350){\blacken\circle{10}}
\put(250,250){\blacken\circle{10}}
\put(250,175){\blacken\circle{10}}
\put(245,150){\llap{\(a_2\)}}
\put(245,225){\llap{\(a_0\)}}
\put(245,325){\llap{\(a\)}}
\put(245,475){\llap{\(a_1\)}}
\Thicklines
\spline(25,175)(50,180)(100,220)(150,350)(200,470)(250,510)(300,470)%
(350,350)(400,220)(450,180)(475,175)
\put(225,25){{\(a\in U_1\)}}
\end{picture}%%%%%%%%%%%%%%%%%%%%%%%%%%%%%%%%%%%%%%%%%%%%%%%%%%%%%
\begin{picture}(625,625)
\put(0,100){\vector(1,0){500}}
\put(250,100){\vector(0,1){450}}
\put(505,110){{\(x\)}}
\put(260,555){{\(g(x)\)}}
\dashline[25]{25}(0,500)(500,500)
\dashline[25]{25}(25,250)(475,250)
\dashline[25]{25}(25,175)(475,175)
\dashline[25]{25}(25,100)(25,500)
\dashline[25]{25}(475,100)(475,500)
\dashline[25]{25}(150,100)(150,250)
\dashline[25]{25}(350,100)(350,250)
\put(25,70){{\(-\pi\)}}\put(475,70){{\(\pi\)}}
\put(150,70){{\(x_1\)}}\put(350,70){{\(x_2\)}}\put(250,70){{\(0\)}}
\put(250,100){\blacken\circle{10}}
\put(250,500){\blacken\circle{10}}
\put(250,350){\blacken\circle{10}}
\put(250,250){\blacken\circle{10}}
\put(250,175){\blacken\circle{10}}
\put(245,150){\llap{\(a_2\)}}
\put(245,225){\llap{\(a\)}}
\put(245,325){\llap{\(a_0'\)}}
\put(245,475){\llap{\(a_1\)}}
\Thicklines
\spline(25,500)(50,490)(100,400)(150,240)(200,190)(250,170)(300,190)%
(350,240)(400,400)(450,490)(475,500)
\put(225,25){{\(a\in U_2\)}}
\end{picture}

\begin{picture}(625,625)(-50,0)
\put(0,100){\vector(1,0){500}}
\put(250,100){\vector(0,1){450}}
\put(505,110){{\(x\)}}
\put(260,555){{\(h(x)\)}}
\dashline[25]{25}(130,300)(370,300)
\dashline[25]{25}(150,100)(150,300)
\dashline[25]{25}(350,100)(350,300)
\path(25,90)(25,110)
\path(475,90)(475,110)
\put(25,70){{\(-\pi\)}}\put(475,70){{\(\pi\)}}
\put(150,70){{\(x_1\)}}\put(350,70){{\(x_2\)}}\put(250,70){{\(0\)}}
\put(250,100){\blacken\circle{10}}
\put(250,450){\blacken\circle{10}}
\put(250,300){\blacken\circle{10}}
\put(245,420){\llap{\(a_2\)}}
\put(245,275){\llap{\(a\)}}
\Thicklines
\spline(90,100)(120,205)(150,310)(200,415)(250,465)(300,415)%
(350,310)(380,205)(410,100)
\put(225,25){{\(a\in U_3\)}}
\end{picture}
\caption{}\label{bil:1}
\end{figure}

В каждой из областей \(U_i\) (\(i=1,2,3\)) введём величины \(\alpha^2\equiv
\alpha^2(a)\), \(\alpha_2^2\equiv\alpha_2^2(a)\) согласно следующим
определяющим их равенствам:
\begin{equation}\label{eq:2:1}
	\begin{aligned}
		\alpha^2&=\dfrac{2}{\pi}\int\limits_{x_1}^{x_2}\sqrt{g(x)-a}\,dx,&
		\alpha_2^2&=\dfrac{2}{\pi}\int\limits_{a>g(x)}\sqrt{a-g(x)}\,dx,&
		&(x\in U_1),\\
		\alpha^2&=\dfrac{2}{\pi}\int\limits_{x_1}^{x_2}\sqrt{a-g(x)}\,dx,&
		\alpha_2^2&=\alpha^2,&&(x\in U_2),\\
		\alpha^2&=\dfrac{2}{\pi}\int\limits_{x_1}^{x_2}\sqrt{h(x)-a}\,dx,&
		\alpha_2^2&=\alpha^2,&&(x\in U_3).
	\end{aligned}
\end{equation}
Заметим, что в области \(U_1\) выполнено \(\alpha^2\to 0\) при \(a\to a_1\),
тогда как в области \(U_2\cup U_3\) выполнено \(\alpha^2\to 0\) при
\(a\to a_2\).

Везде ниже
\begin{equation}\label{eq:2:2}
	b\equiv b(\lambda_1)=\left\{\begin{array}{r}
	\dfrac{1}{2}\lambda\alpha^2\\ -\dfrac{1}{2}\lambda\alpha^2
	\end{array}\right.\qquad
	b_2\equiv b_2(\lambda)=\left\{\begin{array}{l@{\quad}l}
	\dfrac{1}{2}\lambda\alpha_2^2&(a\in U_1\cup U_2)\\
	b&(a\in U_3).\end{array}\right.
\end{equation}

Для фиксированного значения параметра \(a\) через \(\lambda_{n,a}\)
обозначим точки спектра задачи~\eqref{eq:1:1} при их естественном упорядочении
\begin{equation}\label{eq:2:3}
	\ldots\leqslant\lambda_{n-1,a}\leqslant\lambda_{n,a}\leqslant
	\lambda_{n+1,a}\leqslant\ldots\qquad(n\in\mathbb Z^+).
\end{equation}

\begin{theor*}
Пусть \(g\in\mathcal G\) и рассматривается задача~\eqref{eq:1:1}. Тогда
справедливы утверждения~(1), (2).
\begin{enumerate}
\item Положительная часть спектра при \(\lambda\gg 1\), \(\forall a\in U\)
состоит из двух ветвей \(\lambda_{\pm}(a,p)\), нумеруемых выбором знака
\(\pm\) и целым числом \(p\gg 1\), и имеют место равенства
\begin{gather}\label{eq:2:4}
	\lambda_{\pm}(a,p)=\lambda_p^0+F(a)^{-1}H_{\pm}(b_2(\lambda_p^0))
	+R_{\pm}(a,p).\\ \intertext{В этих формулах}\label{eq:2:5}
	\lambda_p^0\equiv\lambda^0(a,p)=2\pi p(F(a))^{-1},\\
	\intertext{величина \(F(a)\) определена в~\eqref{eq:1:8}, имеют
	место оценки}\label{eq:2:6}
	R_{\pm}(a,p)=\left\{\begin{aligned}
	&F(a)^{-1}\,O\left((\lambda_p^0)^{-2/3}\ln\lambda_p^0\right)&
	&(a_2\leqslant a<a_1),\\
	&O\left((\lambda_p^0)^{-1/2}(\ln\lambda_p^0)^{1/2}\right)&
	&(0\leqslant a<a_2),\end{aligned}\right.\\
	\intertext{и функции \(H_{\pm}\) определяются равенствами}
	\label{eq:2:7}
	H_{\pm}(x)=\pm\arctg e^{\pi x}-x+x\ln|x|-
	\arg\Gamma\left(\dfrac12+ix\right),
\end{gather}
в которых \(\Gamma(\cdot)\) "--- гамма-функция, и ветвь аргумента выбрана
таким образом, что \(\arg z=0\;(z\in\mathbb R^+)\) , \(|\arg z|\leqslant
\pi\).

\item Упорядоченные собственные значения \(\lambda_{n,a}\)~\eqref{eq:2:3}
при \(n\gg 1\) удовлетворяют соотношениям
\begin{equation}\label{eq:2:8}
	\begin{aligned}
		\lambda_{n,a}&=\left\{\begin{aligned}
		&\lambda_+\left(a,\dfrac{n}{2}\right)\text{ или }
		\lambda_-\left(a,\dfrac{n}{2}\right)&&(n\text{ "--- чётное})\\
		&\lambda_+\left(a,\dfrac{n+1}{2}\right)\text{ или }
		\lambda_-\left(a,\dfrac{n+1}{2}\right)&&(n\text{ "--- нечётное})
		\end{aligned}\right.&&(a\leqslant a_2-C\lambda^{-1}),\\
		\lambda_{n,a}&=\left\{\begin{aligned}
		&\lambda_+\left(a,\dfrac{n}{2}\right)&&(n\text{ "--- чётное})\\
		&\lambda_-\left(a,\dfrac{n+1}{2}\right)&&(n\text{ "--- нечётное})
		\end{aligned}\right.&&(a\geqslant a_2+C\lambda^{-1}),\\
	\end{aligned}
\end{equation}
и собственные функции, принадлежащие собственным значениям
\(\lambda_{\pm}(a,p)\) имеют \(2p\) нулей на интервале \([-\pi,\pi)\).
\end{enumerate}
\end{theor*}

Свойства функций \(H_{\pm}(x)\) будут рассмотрены ниже. Обратим внимание
на то, что эти функции непрерывны, но их производные имеют
особенность вида \(\ln|x|\) при \(x\to 0\). Функции \(H_{\pm}(x)\)
немонотонны и их качественное поведение ясно из рис.~\ref{bil:2}.

Основой для доказательства сформулированной теоремы является лемма о спектре,
дающая, кроме того, несколько более точную информацию об асимптотике собственных
значений. Для формулировки этой леммы, кроме введённых выше величин \(F(a)\),
\(\alpha^2\equiv\alpha^2(b)\), \(b\equiv b(\lambda,a)\), понадобятся величины
\begin{equation}\label{eq:2:9}
	k(b)=\left(1+e^{2\pi b}\right)^{1/2}-e^{\pi b},
\end{equation}
величины \(\zeta_2\equiv\zeta_2(a)\), определяемые равенствами
\begin{equation}\label{eq:2:10}
	\begin{aligned}
		\int\limits_{x_2}^{\pi}\sqrt{a-g(x)}\,dx&=
		\int\limits_{\alpha}^{\zeta_2}\sqrt{\zeta_2^2-\alpha^2}\,d\zeta=
		\dfrac{1}{2}\zeta_2\sqrt{\zeta_2^2-\alpha^2}-
		\dfrac{\alpha^2}{2}\ln\left(\dfrac{\zeta_2}{\alpha}+
		\dfrac{1}{\alpha}\sqrt{\zeta_2^2-\alpha^2}\right),&
		&(a\in\lambda_1)\\
		\int\limits_{x_2}^{\pi}\sqrt{g(x)-a}\,dx&=
		\int\limits_{\alpha}^{\zeta_2}\sqrt{\zeta^2-\alpha^2}\,d\zeta&
		&(a\in U_2)\\
		\int\limits_0^{\pi}\sqrt{g(x)-a}\,dx&=\int\limits_0^{\zeta_2}
		\sqrt{\zeta^2+\alpha^2}\,d\zeta=\dfrac{1}{2}\zeta_2
		\sqrt{\zeta_2^2+\alpha^2}+\dfrac{\alpha^2}{2}\ln
		\left(\dfrac{\zeta_2}{\alpha}+\dfrac{1}{\alpha}
		\sqrt{\zeta_2^2+\alpha^2}\right)&
		&(a\in U_3),
	\end{aligned}
\end{equation}
и величина
\begin{equation}\label{eq:2:11}
	\Psi(\lambda,a)=\lambda\zeta_2^2-2b\ln(\zeta_2\sqrt{2\lambda})+
	\arg\Gamma\left(\dfrac{1}{2}+ib\right).
\end{equation}
Введём более детальное, чем \(U=\bigcup_i U_i\), разбиение интервала
\(U\) "--- \(U=\bigcup_{i=1}^5 A_i\). По определению
\begin{align*}
	a&\in A_1=U_1&&(a_0\leqslant a<a_1),\\
	a&\in A_2\subset U_2&&(a_2+C\lambda^{-1+\varepsilon_1}\leqslant
	a\leqslant a_0',\quad 0<\varepsilon_1<1,\quad C>0),\\
	a&\in A_3\subset U_2&&(a_2\leqslant a\leqslant a_2+C\lambda^{-1+
	\varepsilon_2},\quad 0<\varepsilon_2<1/2,\quad C>0),\\
	a&\in A_4\subset U_3&&(a_2-C\lambda^{-1+\varepsilon_2'}\leqslant
	a\leqslant a_2,\quad 0<\varepsilon_2'<1/2,\quad C>0),\\
	a&\in A_5\subset U_3&&(0\leqslant a\leqslant a_2-C\lambda^{-1+
	\varepsilon_1'},\quad 0<\varepsilon_1'<1,\quad C>0).
\end{align*}

\begin{los}
Пусть \(g\in\mathcal G\), \(a\in U\) и \(\lambda>C\gg 1\) "--- любое
положительное собственное значение для задачи~\eqref{eq:1:1}. Тогда для
любого такого \(\lambda\) существуют целое число \(p\in\mathbb Z^+\)
(\(p>C\gg 1\)) и выбор знака \(\pm\) такие, что выполняются
утверждения~(1)--(5):
\begin{enumerate}
\item Если \(a\in A_1\), то имеет место равенство
\begin{equation}\label{eq:2:12}
	\lambda F(a)=2\pi p\pm\dfrac{\pi}{2}+
	O\left(\lambda^{-2/3}\ln\lambda\right).
\end{equation}
\item Если \(a\in A_2\), то \(\forall\varepsilon_1\) такого, что
\(0<\varepsilon_1<1\), имеет место равенство
\begin{equation}\label{eq:2:13}
	\lambda F(a)=2\pi p\pm\dfrac{\pi}{2}-\dfrac{1}{24 b}+
	O\left(b^{-3}+\lambda^{-2/3}\ln\lambda\right).
\end{equation}
\item Если \(a\in A_3\), то \(\forall\varepsilon_2\) такого, что
\(0<\varepsilon_2<1/2\), имеет место равенство
\begin{equation}\label{eq:2:14}
	\Psi(\lambda,a)=2\pi p\pm\arccos\dfrac{2k(b)}{1+k^2(b)}-
	\dfrac{b^2}{2\zeta_2^2\lambda}+O\left(\dfrac{b^4}{\lambda^2}+
	\lambda^{-2/3}\ln\lambda\right).
\end{equation}
\item Если \(a\in A_4\), то \(\forall\varepsilon_2'\) такого, что
\(0<\varepsilon_2'<1/2\), имеет место равенство
\begin{equation}\label{eq:2:15}
	\Psi(\lambda,a)=2\pi p\pm\arccos\dfrac{2k(b)}{1+k^2(b)}-
	\dfrac{b^2}{2\zeta_2^2\lambda}+O\left(\dfrac{b^4}{\lambda^2}+
	\lambda^{-1}\ln\lambda\right).
\end{equation}
\item Если \(a\in A_5\), то \(\forall\varepsilon_1'\) такого, что
\(0<\varepsilon_1'<1\), имеет место равенство
\begin{equation}\label{eq:2:16}
	\lambda F(a)=2\pi p-\dfrac{1}{24 b}+O\left(|b|^{-3/2}+
	\lambda^{-1/2}\,(\ln\lambda)^{1/2}\right).
\end{equation}
\end{enumerate}
\end{los}

Рассмотрим свойства функций \(H_{\pm}(x)\).

В случае \(|x|\gg 1\) воспользуемся формулой Стирлинга \cite{22}
\begin{gather*}
	\ln\Gamma(z)=\left(z-\dfrac{1}{2}\right)\ln z-z+\dfrac{1}{2}
	\ln 2\pi+\sum\limits_{n=1}^m\dfrac{B_{2n}}{2n(2n-1)z^{2n-1}}+
	O\left(|z|^{-2m-1}\right),\\
	B_2=1/6,\qquad B_4=-1/30,\qquad (|\arg z|<\pi).
\end{gather*}
Так как \(\arg\Gamma(1/2+ix)=\Im\ln\Gamma(1/2+ix)\), то
\[
	\arg\Gamma(1/2+ix)=x\ln|x|-x+\dfrac{1}{24x}+
	\dfrac{7}{24\cdot120}\dfrac{1}{x^3}+O(|x|^{-5})
\]
и следовательно при \(|x|\gg 1\)
\begin{equation}\label{eq:2:17}
	\begin{aligned}
		H_{\pm}(x)&=\pm\dfrac{\pi}{2}-\dfrac{1}{24x}-\dfrac{7}{24\cdot
		120}\dfrac{1}{x^3}+O\left(|x|^{-5}\right)&&(x>0),\\
		H_{\pm}(x)&=\dfrac{1}{24x}-\dfrac{7}{24\cdot 120}\dfrac{1}{x^3}+
		O\left(|x|^{-5}\right)&&(x<0).
	\end{aligned}
\end{equation}
В случае \(|x|\ll 1\) можно воспользоваться следующей формулой (см.~\cite{23})
\[
	\arg\Gamma(x+iy)=y\psi\left(\dfrac{1}{2}\right)+\sum\limits_{n=0}^{\infty}
	\left(\dfrac{y}{x+n}-\arctg\dfrac{x}{2n+1}\right)\qquad
	\left(\psi(x)=\dfrac{\Gamma'(x)}{\Gamma(x)}\right),
\]
из которой следует, что
\[
	\arg\Gamma\left(\dfrac{1}{2}+ix\right)=x\psi\left(\dfrac{1}{2}\right)+
	\sum\limits_{n=0}^{\infty}\left(\dfrac{2x}{2n+1}-\arctg\dfrac{2x}{2n+1}
	\right).
\]
Разлагая функцию \(\arctg\dfrac{2x}{2n+1}\) в ряд Тейлора, получим
\begin{gather*}
	\arg\Gamma\left(\dfrac{1}{2}+ix\right)=x\psi\left(\dfrac{1}{2}\right)+
	\sum\limits_{n=1}^{\infty}A_{2n+1}\cdot x^{2n+1},\\
	A_{2n+1}=(-1)^{n+1}\dfrac{2^{2n+1}}{2n+1}\sum\limits_{m=0}^{\infty}
	(2m+1)^{-2n-1}\qquad \left(\psi\left(\dfrac{1}{2}\right)=\gamma-
	2\ln 2=-1,9635\ldots\right).
\end{gather*}
Так как
\[
	\arctg e^{\pi x}=\dfrac{\pi}{4}-\dfrac{\pi x}{2}+\dfrac{1}{2}(\pi x)^3+
	O\left(|x|^5\right),
\]
то при \(|x|\ll 1\)
\begin{equation}\label{eq:2:18}
	H_{\pm}(x)=\pm\dfrac{\pi}{4}+x\ln|x|+x\left(\pm\dfrac{\pi}{2}-1-
	\psi\left(\dfrac{1}{2}\right)\right)+x^3\left(\mp\dfrac{\pi^3}{12}-
	A_3\right)+O\left(|x|^5\right).
\end{equation}
Прямо из определений~\eqref{eq:2:7} функций \(H_{\pm}(x)\) следует, что
\begin{equation}\label{eq:2:19}
	H_{+}(x)=H_{-}(x)+2\arctg e^{\pi x},\qquad H_{\pm}(-x)=\pm\dfrac{\pi}{2}-
	H_{\pm}(x)
\end{equation}
и таким образом достаточно проанализировать поведение функции \(H_+(x)\)
при \(x\geqslant 0\). Из формул~\eqref{eq:2:17},~\eqref{eq:2:18} следует, что
при \(x\geqslant 0\) эта функция имеет минимум в точке \(x_{min}^+\simeq
0,0293\) и \(H_+(x_{\min}^+)\simeq\pi/4-x_{min}^+\). Функция \(H_{-}(x)\)
также имеет один минимум в точке \(x_{min}^-\simeq 1,683\) и
\(H_-(x_{min}^-)\simeq\pi/2-0,0203\). Вид функций \(H_{\pm}(x)\)
представлен на рис.~\ref{bil:2}.

\begin{figure}[ht]
\unitlength=0.132mm
\begin{picture}(1250,1250)
\put(425,50){\vector(0,1){1150}}
\put(50,625){\vector(1,0){1150}}
\put(425,625){\blacken\circle{10}}
\put(425,825){\blacken\circle{10}}
\put(425,1025){\blacken\circle{10}}
\put(425,425){\blacken\circle{10}}
\put(425,225){\blacken\circle{10}}
\put(475,625){\blacken\circle{10}}
\put(860,625){\blacken\circle{10}}
\dashline[25]{25}(860,625)(860,175)
\dashline[25]{25}(425,225)(1150,225)
\dashline[25]{25}(425,1025)(1150,1025)
\dashline[25]{25}(425,175)(900,175)
\dashline[25]{25}(475,625)(475,790)
\dashline[25]{25}(425,790)(825,790)
\put(420,600){\llap{\(0\)}}
\put(420,1020){\llap{\(\dfrac{\pi}{2}\)}}
\put(420,220){\llap{\(-\dfrac{\pi}{2}\)}}
\put(420,820){\llap{\(\dfrac{\pi}{4}\)}}
\put(420,420){\llap{\(-\dfrac{\pi}{4}\)}}
\put(480,600){{\(x^+_{min}\simeq 0,03\)}}
\put(875,600){{\(x^-_{min}\simeq 1,68\)}}
\put(835,785){{\(H_+(x^+_{min})\simeq\dfrac{\pi}{4}-0,03\)}}
\put(915,165){{\(H_-(x^-_{min})\simeq-\dfrac{\pi}{2}-0,02\)}}
\Thicklines
\spline(50,630)(120,650)(250,700)(325,635)(410,500)(420,470)(425,425)(430,380)%
(440,350)(500,275)(600,225)(650,212)(700,200)(750,190)(800,180)(860,175)%
(920,190)(1000,215)(1150,220)
\put(600,250){{\(H_-(x)\)}}
\spline(50,640)(120,660)(250,740)(300,800)(350,860)(375,885)(400,880)%
(415,861)(417,849)(420,840)(425,825)(427,810)(450,795)(460,792)(475,790)%
(525,794)(600,920)(800,985)(900,1005)(1150,1022)
\put(775,935){{\(H_+(x)\)}}
\end{picture}
\caption{}\label{bil:2}
\end{figure}

С качественной точки зрения полученный результат означает, что
в рассматриваемой задаче величины \(s_{n,a}=F(a)\lambda_{n,a}\) при \(n\gg 1\),
\(a\in U\) ведут себя так же, как величины \(\sqrt{\tilde\lambda_n}\)
в классической (дефинитной) задаче Штурма--Лиувилля, т.~е. группируются
в пары, и между этими парами имеются лакуны, ширина \(\Delta\) которых
в рассматриваемой задаче имеет асимптотику вида
\begin{equation}\label{eq:2:20}
	F(a)\cdot\Delta=2\pi-2\arctg e^{\pi b_2(\lambda)}+O\left(\lambda^{-1/2}
	(\ln\lambda)^{1/2}\right)\geqslant\pi+H_{+}(x_{min}^-)+
	O\left(\lambda^{-1/2}(\ln\lambda)^{1/2}\right).
\end{equation}

\section{Результаты Олвера}\label{pt:3}
В этом параграфе будут кратко описаны результаты работы \cite{17} в той
их части, которая нам необходима. Попутно будут введены нужные определения,
обозначения, и будет доказано, что результаты работы \cite{17} применимы
при \(g\in\mathcal G\). Как уже отмечалось, класс \(\mathcal G\) и был
определён таким образом, чтобы для него были верны результаты работы
\cite{17}.

Точки поворота \(x_i\equiv x_i(a)\), \(x_1=-x_2\), \(x_2>0\) были определены
выше, и при \(g\in\mathcal G\) функция \(x_2(a)\) имеет обратную \(a=a(x_2)\).
Рассмотрим функцию
\begin{equation}\label{eq:3:1}
	f(x,a)=a-g(x).
\end{equation}
Для применимости результатов Олвера достаточно, чтобы имели место
представления
\begin{align*}
	f(x,a)&=(x^2-x_2^2)\rho_1(x,x_2)&&(a\in U_1),\\
	f(x,a)&=(x_2^2-x^2)\rho_2(x,x_2)&&(a\in U_2),\\
	f(x,a)&=-(x^2+x_2^2)\rho_3(x,x_2)&&(a\in U_3),
\end{align*}
в которых функции \(\rho_i(x,x_2)\geqslant C>0\) и имеют по крайней мере
четыре ограниченные производные по \(x\). Кроме того, предполагается, что
функция \(f(x,a(x_2))\) не возрастает с ростом \(x_2\) при \(a\in U_1\cup
U_3\) и неубывает "--- при \(a\in U_2\). Докажем, что эти условия выполняются
при \(g\in\mathcal G\) и \(a\in U\). Пусть \(a\in U_1\). Так как \(x_2\)
убывает с ростом \(a\) (см.~рис.~\ref{bil:1}), то прямо из
определения~\eqref{eq:3:1} следует, что функция \(f(x,a(x_2))\) убывает
с ростом \(x_2\). Представление \(f(x,a)=(x^2-x_2^2)\rho_1(x,x_2)\)
следует из формулы Тейлора
\[
	f(x,a)=(x-x_2)\int\limits_0^1f^{(1)}(x_2+s(x-x_2))\,ds=(x-x_2)f_1(x,a).
\]
Применяя затем эту формулу с заменой \(x_2\) на \(x_1=-x_2\) к функции
\(f_1(x,a)\), получим, что при \(a\in U_2\)
\[
	\rho_1(x,x_2)=\int\limits_0^1\,dt\int\limits_0^1\,ds
	f^{(2)}(x_2+s(-2x_2+t(x+x_2)).
\]
Условие \(\rho_1(x,x_2)\geqslant C>0\) следует из простоты нулей \(x_i\)
функции \(f\) и того, что \(g^{(2)}(0)\leqslant -C\) (\(C>0\)). Таким образом,
при \(a\in U_1\) применимы результаты Олвера. При \(a\in U_2\) доказательство
аналогично и надо просто заменить \(f\) на \(-f\). При \(a\in U_3\) эти же
рассуждения надо применить к функции \(\tilde f(x,a)=a-h(x)\). Таким образом,
при \(g\in\mathcal G\) и любом \(a\in U\) можно применять результаты работы
\cite{17}.

Метод Олвера построения асимптотик фундаментальной системы решений
уравнения~\eqref{eq:1:1} основан на применении преобразования Лиувилля, более
сложного, чем указанное в~\ref{pt:1}. Нужное преобразование (см.
\cite{1,2}) определяется соотношениями
\begin{equation}\label{eq:3:2}
	u(x)=(\varphi'(\zeta))^{1/2}\,w(\zeta),\qquad x=\varphi(\zeta),
\end{equation}
где \(w\) "--- новая зависимая, а \(\zeta\) "--- новая независимая
переменные. В новых переменных уравнение~\eqref{eq:1:1} принимает вид
\begin{equation}\label{eq:3:3}
	w''(\zeta)=\left(\lambda^2f(\varphi')^2-\dfrac{1}{2}\{\varphi,
	\zeta\}\right)\,w(\zeta),
\end{equation}
где \(\{\varphi,\zeta\}\) "--- производная Шварца (\(\{\varphi,\zeta\}=
\varphi'''(\varphi')^{-1}-3/2\cdot(\varphi'')^2(\varphi')^{-2}\)).
В этих формулах и везде ниже штрих означает производную по \(\zeta\).

Рассмотрим случай \(a\in U_1\). Тогда функция \(\varphi\) определяется
из условий
\begin{equation}\label{eq:3:4}
	(\varphi')^2\,f=\zeta^2-\alpha^2,\qquad\varphi(\alpha)=x_2,\qquad
	(\varphi(-\alpha)=x_1=-x_2),
\end{equation}
в которых величина \(\alpha=\alpha(a)\) определена в~\eqref{eq:2:1}. Отсюда
получаем, что переменные \(\zeta\) и \(x\) при \(a\in U_1\) связаны
соотношениями
\begin{equation}\label{eq:3:5}
	\begin{aligned}
		\int\limits_{x_1}^x\sqrt{-f}\,dx&=\int\limits_{-\alpha}^{\zeta}
		\sqrt{\alpha^2-\zeta^2}\,d\zeta&&(x_1<x<x_2),\\
		\int\limits_{x}^{x_1}\sqrt{f}\,dx&=\int\limits_{\zeta}^{-\alpha}
		\sqrt{\zeta^2-\alpha^2}\,d\zeta&&(x<x_1),\\
		\int\limits_{x_2}^x\sqrt{f}\,dx&=\int\limits_{\alpha}^{\zeta}
		\sqrt{\zeta^2-\alpha^2}\,d\zeta&&(x>x_2).
	\end{aligned}
\end{equation}
В области \(U_2\) функция \(\varphi(\zeta)\) определяется из уравнения
\begin{equation}\label{eq:3:6}
	(\varphi')^2\,f=\alpha^2-\zeta^2,
\end{equation}
и всё аналогично случаю \(a\in U_1\), надо только заменить в~\eqref{eq:3:5}
\(f\) на \(-f\). В области \(U_3\) функция \(\varphi\) определяется
из условий
\begin{equation}\label{eq:3:7}
	(\varphi')^2\,f=-\alpha^2-\zeta^2,\qquad\varphi(0)=0,
\end{equation}
и, следовательно,
\begin{equation}\label{eq:3:8}
	\int\limits_0^x\sqrt{-f}\,dx=\int\limits_0^{\zeta}\sqrt{\alpha^2+
	\zeta^2}\,d\zeta.
\end{equation}
При всех \(a\in U\) функция \(x=\varphi(\zeta)\) имеет обратную функцию
\(\zeta=\zeta(x)\), которая монотонно возрастает на отрезке \(|x|\leqslant
\pi\) и
\begin{equation}\label{eq:3:9}
	\zeta_2=\zeta(\pi)=-\zeta_1=\zeta(-\pi),\qquad\zeta_2>C>0.
\end{equation}
Отбрасывая в уравнении~\eqref{eq:3:3} член, содержащий \(\{\varphi,\zeta\}\),
получим в соответствии с равенствами~\eqref{eq:3:4}, \eqref{eq:3:6}
и~\eqref{eq:3:7} следующие модельные уравнения:
\begin{equation}\label{eq:3:10}
	\begin{aligned}
		v''(\zeta)&=\lambda^2(\zeta^2-\alpha^2)v(\zeta)&&(a\in U_1),\\
		v''(\zeta)&=\lambda^2(\alpha^2-\zeta^2)v(\zeta)&&(a\in U_2),\\
		v''(\zeta)&=-\lambda^2(\alpha^2+\zeta^2)v(\zeta)&&(a\in U_3).		
	\end{aligned}
\end{equation}
Линейной заменой координат эти уравнения приводятся к уравнениям Вебера
\begin{equation}\label{eq:3:11}
	\dfrac{d^2y}{dz^2}(z)=\left(\beta\pm\dfrac{z^2}{4}\right)y(z).
\end{equation}
Различные частные решения этого уравнения называются функциями Вебера
или функциями параболического цилиндра. Ниже, следуя работе \cite{17},
будут использоваться две конкретные системы линейно независимых решений
уравнения~\eqref{eq:3:11} "--- \(U(\beta,z)\), \(\overline{U}(\beta,z)\)
и \(W(\beta,\pm z)\). Все нужные нам результаты о свойствах этих решений
содержатся в работах \cite{17,24,25}. Функции \(U(\beta,z)\),
\(\overline{U}(\beta,z)\) и \(W(\beta,\pm z)\) следующим образом связаны
со стандартными функциями параболического цилиндра \(D_{\nu}(z)\)
и функциями Уиттекера \(W_{\varkappa,\nu}(z)\) \cite{26}:
\begin{align*}
	U(\beta,z)&=D_{-\beta-1/2}(z)=2^{-\beta/2}z^{-1/2}
	W_{-1/2\beta,-1/4}(z),\\
	\overline{U}(\beta,z)&=\tg\pi\beta\cdot U(\beta,z)-\sec\pi\beta
	U(\beta,-z),\\
	W(\beta,z)&=\sqrt{2k(\beta)}e^{i\pi\beta/4}\Re\left\{
	e^{i\left(\phi/2+1/8\right)}U\left(i\beta,e^{-i\pi/4}z\right)\right\},\\
	W(\beta,-z)&=\sqrt{\dfrac{2}{k(\beta)}}e^{i\pi\beta/4}\Im\left\{
	e^{i\left(\phi/2+1/8\right)}U\left(i\beta,e^{-i\pi/4}z\right)\right\}.
\end{align*}
В этих формулах \(\phi=\arg\Gamma\left(1/2+i\beta\right)\), и величина
\(k(\beta)\) определена в~\eqref{eq:2:9}. При \(a\in U_1\) первое
из уравнений~\eqref{eq:3:10} приводится к уравнению~\eqref{eq:3:11} со знаком
плюс в правой части, и в качестве фундаментальной системы решений
уравнения~\eqref{eq:3:10} для \(v(\zeta)\) удобно выбрать
\begin{align*}
	v_1(\zeta)&=U(\beta,z),& v_2(\zeta)&=\overline{U}(\beta,z),\\
	\beta&=-\dfrac{1}{2}\lambda\alpha^2,& z&=\zeta\sqrt{2\lambda}.
\end{align*}
В качестве фундаментальной системы решений можно выбрать и функции \(U(\beta,
-z)\), \(\overline{U}(\beta,-z)\), и между этими системами имеются уравнения
связи
\begin{equation}\label{eq:3:12}
	\begin{aligned}
		U(\beta,-z)&=\cos(\pi\beta)\overline{U}(\beta,-z)-
		\sin(\pi\beta)U(\beta,z),\\
		\overline{U}(\beta,-z)&=\sin(\pi\beta)\overline{U}(\beta,z)+
		\cos(\pi\beta)U(\beta,z).
	\end{aligned}
\end{equation}
В случае \(a\in U_2\cup U_3\) соответствующее уравнение~\eqref{eq:3:10}
приводится к уравнению~\eqref{eq:3:11} со знаком минус в правой части,
и в качестве фундаментальных систем решений уравнений~\eqref{eq:3:10}
выбираются, соответственно,
\begin{align*}
	v_i(\zeta)&=W\left(\dfrac{1}{2}\lambda\alpha^2\pm\zeta\sqrt{2\lambda}
	\right)&i&=1,2,&&(a\in U_2),\\
	v_i(\zeta)&=W\left(-\dfrac{1}{2}\lambda\alpha^2\pm\zeta\sqrt{2\lambda}
	\right)&i&=1,2,&&(a\in U_3).
\end{align*}
Из результатов работы \cite{17} вытекают следующие нужные нам утверждения.

Пусть \(a\in U_1\subset [0,\zeta_2]\). Тогда уравнение~\eqref{eq:3:3} при
\(g\in\mathcal G\) имеет два линейно независимых решения вида
\begin{equation}\label{eq:3:13}
	\begin{aligned}
		w_1(\zeta)&=U\left(-\dfrac{1}{2}\lambda\alpha^2,\zeta
		\sqrt{2\lambda}\right)+R_1(\zeta),\\
		w_2(\zeta)&=U\left(-\dfrac{1}{2}\lambda\alpha^2,\zeta
		\sqrt{2\lambda}\right)+R_2(\zeta),
	\end{aligned}
\end{equation}
и при \(\zeta\in[0,\zeta_2]\) для функций \(R_i(\zeta)\) и их производных
\(R_i'(\zeta)\) справедливы следующие оценки
\begin{equation}\label{eq:3:14}
	\begin{aligned}
		R_1(\zeta)&=\left(E(\beta,z)\right)^{-1}M(\beta,z)
		O\left(\lambda^{-2/3}\ln\lambda\right),\\
		R_2(\zeta)&=E(\beta,z)M(\beta,z)
		O\left(\lambda^{-2/3}\ln\lambda\right),\\
		R_1'(\zeta)&=\left(E(\beta,z)\right)^{-1}N(\beta,z)
		O\left(\lambda^{-1/6}\ln\lambda\right),\\
		R_2'(\zeta)&=E(\beta,z)N(\beta,z)
		O\left(\lambda^{-1/6}\ln\lambda\right).
	\end{aligned}
\end{equation}
В этих равенствах \(\beta=-\dfrac{1}{2}\lambda\alpha^2\), \(z=\zeta
\sqrt{2\lambda}\). При \(a\in U\) нам понадобятся только значения функций
\(E(\beta,z)\), \(M(\beta,z)\) и \(N(\beta,z)\) при \(\beta>0\) и
\(z>\rho(\beta)\), где \(\rho(\beta)\) "--- максимальный действительный
корень уравнения \(U(\beta,z)=\overline{U}(\beta,z)\), и при
\(|\beta|\gg 1\), \(\beta<0\), выполняется соотношение
\begin{equation}\label{eq:3:15}
	\rho(\beta)=2\sqrt{-\beta}+C(-\beta)^{1/2}+O\left(|\beta|^{-5/6}
	\right)\qquad (C\simeq-0,366\ldots).
\end{equation}
Если \(a\in U_1\), \(\beta<0\), \(z>\rho(\beta)\), то
в равенствах~\eqref{eq:3:14}
\begin{equation}\label{eq:3:16}
	\begin{gathered}
		E(\beta,z)=\left(\overline{U}(\beta,z)\right)^{1/2}
		\left(U(\beta,z)\right)^{-1/2},\qquad
		M(\beta,z)=\left[2U(\beta,z)\overline{U}(\beta,z)\right]^{1/2},\\
		N(\beta,z)=\left[\left(U^{(1)}(\beta,z)\right)^2
		\dfrac{\overline{U}(\beta,z)}{U(\beta,z)}+
		\left(\overline{U}^{(1)}(\beta,z)\right)^2
		\dfrac{U(\beta,z)}{\overline{U}(\beta,z)}\right]^{1/2}.
	\end{gathered}
\end{equation}
Здесь и ниже \(U^{(1)}(\beta,z)=\dfrac{d}{dz}U(\beta,z)\),
\(\overline{U}^{(1)}(\beta,z)=\dfrac{d}{dz}\overline{U}(\beta,z)\).
Что касается оценки остаточных членов при \(\zeta<0\), то, как показано
в работе \cite{17}, существуют два линейно независимых решения \(w_3(\zeta)\),
\(w_4(\zeta)\) уравнения~\eqref{eq:3:3} такие, что при \(a\in U_1\)
\begin{equation}\label{eq:3:17}
	\begin{aligned}
		w_3(\zeta)&=U\left(-\dfrac{1}{2}\lambda\alpha^2,
		-\zeta\sqrt{2\lambda}\right)+R_3(\zeta),\\
		w_4(\zeta)&=\overline{U}\left(-\dfrac{1}{2}\lambda\alpha^2,
		-\zeta\sqrt{2\lambda}\right)+R_4(\zeta),
	\end{aligned}
\end{equation}
и величины \(R_3(\zeta)\), \(R_3'(\zeta)\) при \(\zeta\in[\zeta_1,0]\)
удовлетворяют тем же равенствам~\eqref{eq:3:14}, что и \(R_1(\zeta)\),
\(R_1'(\zeta)\), а \(R_4(\zeta)\), \(R_4'(\zeta)\) "--- тем же
равенствам~\eqref{eq:3:14}, что и \(R_2(\zeta)\), \(R_2'(\zeta)\), но в
равенствах~\eqref{eq:3:14} теперь надо положить \(z=-\zeta\sqrt{2\lambda}\).

При \(a\in U_1\) для коэффициентов \(A_1\), \(B_1\), \(A_2\), \(B_2\)
в уравнениях связи
\begin{equation}\label{eq:3:18}
	\begin{aligned}
		w_1(\zeta)&=A_1w_3(\zeta)+B_1w_4(\zeta),\\
		w_2(\zeta)&=A_2w_3(\zeta)+B_2w_4(\zeta),
	\end{aligned}
\end{equation}
имеют место равенства
\begin{equation}\label{eq:3:19}
	\begin{aligned}
		A_1&=\sin\left(\dfrac{\pi}{2}\lambda\alpha^2\right)+
		O\left(\lambda^{-2/3}\ln\lambda\right),&
		B_1&=\cos\left(\dfrac{\pi}{2}\lambda\alpha^2\right)+
		O\left(\lambda^{-2/3}\ln\lambda\right),\\
		A_2&=\cos\left(\dfrac{\pi}{2}\lambda\alpha^2\right)+
		O\left(\lambda^{-2/3}\ln\lambda\right),&
		B_2&=-\sin\left(\dfrac{\pi}{2}\lambda\alpha^2\right)+
		O\left(\lambda^{-2/3}\ln\lambda\right).
	\end{aligned}
\end{equation}
Переходим к рассмотрению случая \(a\in U_2\). В этом случае нужная нам часть
результатов Олвера формулируется следующим образом: при \(g\in\mathcal G\),
\(a\in U_2\) имеются два линейно независимых решения \(w_1\), \(w_2\)
уравнения~\eqref{eq:3:3} такие, что
\begin{equation}\label{eq:3:20}
	\begin{aligned}
		w_1(\zeta)&=\left(k(b)\right)^{-1/2}W(b,\zeta\sqrt{2\lambda})+
		R_1(\zeta),&&\left(b=\dfrac{1}{2}\lambda\alpha^2\right),\\
		w_2(\zeta)&=\left(k(b)\right)^{1/2}W(b,-\zeta\sqrt{2\lambda})+
		R_2(\zeta),
	\end{aligned}
\end{equation}
и при этом для \(\zeta\in[0,\zeta_2]\) снова выполняются
равенства~\eqref{eq:3:14}, в которых \(\beta=b=\dfrac{1}{2}\lambda\alpha^2\),
\(z=\zeta\sqrt{2\lambda}\), \(E(\beta,z)=1\). Нам будут нужны только выражения
\(M(\beta,z)\), \(N(\beta,z)\) при \(z\geqslant\sigma(\beta)\), где
\(\sigma(\beta)\) "--- минимальный положительный корень уравнения \(k(b)^{-1/2}
W(b,-z)=k(b)^{1/2}W(b,-z)\), и при \(b\gg 1\)
\begin{equation}\label{eq:3:21}
	\sigma(b)=2b^{1/2}-cb^{-1/6}+O\left(b^{-5/6}\right)\qquad
	(C\simeq-0,366\ldots)
\end{equation}
Если \(z\geqslant\sigma(b)\), \(\beta=\beta\geqslant 0\) и \(a\in U_2\),
то в равенствах~\eqref{eq:3:14}
\begin{equation}\label{eq:3:22}
	\begin{aligned}
		M(\beta,z)&=\left[k(b)^{-1}W^2(b,z)+k(b)W^2(b,-z)\right]^{-1/2},\\
		N(\beta,z)&=\left[k(b)^{-1}\left(W^{(1)}(b,z)\right)^2+
		k(b)\left(W^{(1)}(b,-z)\right)^2\right]^{-1/2},\\
	\end{aligned}
\end{equation}
Кроме того, при \(a\in U_2\) имеются два линейно независимых решения \(w_3\),
\(w_4\) уравнения~\eqref{eq:3:3}
\begin{equation}\label{eq:3:23}
	\begin{aligned}
		w_3(\zeta)&=\left(k(b)\right)^{-1/2}W(b,-\zeta\sqrt{2\lambda})+
		R_3(\zeta),\\
		w_4(\zeta)&=\left(k(b)\right)^{1/2}W(b,\zeta\sqrt{2\lambda})+
		R_4(\zeta),
	\end{aligned}
\end{equation}
для которых при \(\zeta\in[\zeta_1,0]\) величины \(R_{3,4}(\zeta)\),
\(R_{3,4}'(\zeta)\) оцениваются точно так же, как описано выше при \(a\in U_1\).
В уравнениях связи~\eqref{eq:3:18} при \(a\in U_2\)
\begin{equation}\label{eq:3:24}
	\begin{aligned}
		A_1&=O\left(\lambda^{-2/3}\ln\lambda\right),&
		B_1&=k(b)^{-1}\left(1+O\left(\lambda^{-2/3}\ln\lambda\right)\right),\\
		A_2&=k(b)\left(1+O\left(\lambda^{-2/3}\ln\lambda\right)\right),&
		B_2&=O\left(\lambda^{-2/3}\ln\lambda\right).
	\end{aligned}
\end{equation}
Наконец, при \(a\in U_3\) в работе \cite{17} доказано, что на всём отрезке
\(\zeta\in[\zeta_1,\zeta_2]\) имеется два линейно независимых решения \(w_1\),
\(w_2\) уравнения~\eqref{eq:3:3}
\begin{equation}\label{eq:3:25}
	\begin{aligned}
		w_1(\zeta)&=W(b,\zeta\sqrt{2\lambda})+R_1(\zeta),&
		\left(b=-\dfrac{1}{2}\lambda\alpha^2\right),\\
		w_2(\zeta)&=W(b,-\zeta\sqrt{2\lambda})+R_2(\zeta)
	\end{aligned}
\end{equation}
с оценками остаточных членов
\begin{equation}\label{eq:3:26}
	R_{1,2}(\zeta)=M(b,z)O\left(\lambda^{-1}\ln\lambda\right),\qquad
	R_{1,2}'(\zeta)=N(b,z)O\left(\lambda^{-1}\ln\lambda\right).
\end{equation}
Величины \(M(b,z)\) и \(N(b,z)\) при этом определяются
соотношениями~\eqref{eq:3:12}, в которых \(\beta=b=-\dfrac{1}{2}\lambda
\alpha^2\), \(z=\zeta\sqrt{2\lambda}\).

\section{Уравнение для спектра}\label{pt:4}
Уравнение для спектра получается из условия существования нетривиальных
периодических решений уравнения~\eqref{eq:1:1}, для которых должны выполняться
условия
\begin{equation}\label{eq:4:1}
	u(-\pi)=u(\pi),\qquad\dfrac{du}{dx}(-\pi)=\dfrac{du}{dx}(\pi).
\end{equation}
Переходя к переменным \(\zeta\), \(w\) \eqref{eq:3:2} и учитывая
определения~\eqref{eq:3:9} величин \(\zeta_i\) (\(i=1,2\)), перепишем
условия~\eqref{eq:4:1} в виде
\begin{gather}\label{eq:4:2}
	w(\zeta_1)\left(\varphi'(\zeta_1)\right)^{1/2}=
	w(\zeta_2)\left(\varphi'(\zeta_2)\right)^{1/2},\\ \label{eq:4:3}
	\dfrac{1}{2}\dfrac{\varphi''(\zeta_1)}{\left(
	\varphi'(\zeta_1)\right)^{3/2}}w(\zeta_1)+\dfrac{1}{\left(
	\varphi'(\zeta_1)\right)^{1/2}}w'(\zeta_1)=
	\dfrac{1}{2}\dfrac{\varphi''(\zeta_2)}{\left(
	\varphi'(\zeta_2)\right)^{3/2}}w(\zeta_2)+\dfrac{1}{\left(
	\varphi'(\zeta_2)\right)^{1/2}}w'(\zeta_2).
\end{gather}
Из периодичности функции \(f\) \eqref{eq:3:1} и равенств~\eqref{eq:3:9}
следует, что во всей области \(U\) выполняется равенство
\(\varphi'(\zeta_1)=\varphi'(\zeta_2)\). Так как \(\varphi'(\zeta)>C>0\),
то уравнение~\eqref{eq:4:2} приводит к равенству \(w(\zeta_1)=w(\zeta_2)\).
Пусть \(a\in U_1\). Дифференцируя уравнение~\eqref{eq:3:4} и учитывая, что
в силу чётности функции \(g\in G\) выполняются равенства \(\dfrac{df}{dx}
(\pm\pi)=0\), получим равенство \(\varphi''(\zeta_2)=\zeta_2\left(
f(\pi)\varphi'(\zeta_2)\right)^{-1}=-\varphi''(\zeta_1)\). Так как
\(\zeta_1=-\zeta_2\), то в области \(U_1\) условия~\eqref{eq:4:2},
\eqref{eq:4:3} принимают вид
\begin{equation}\label{eq:4:4}
	\begin{aligned}
		w(-\zeta_2)&=w(\zeta_2)\\
		w'(-\zeta_2)&=w'(\zeta_2)+\dfrac{\zeta_2}{\zeta_2^2-\alpha^2}
		w(\zeta_2).
	\end{aligned}
\end{equation}
При выводе второго из этих условий было использовано равенство
\[
	\varphi''(\zeta_2)\left(\varphi'(\zeta_2)\right)^{-1}=
	\zeta_2\left(\varphi'(\zeta_2)\right)^2f(\pi)^{-1}=\zeta_2
	\left(\zeta_2^2-\alpha^2\right)^{-1}.
\]
Точно так же получим, что условия~\eqref{eq:4:4} имеют место и в области
\(U_2\), а в области \(U_3\) надо только в правой части второго из условий
заменить \(\zeta_2^2-\alpha^2\) на \(\zeta_2^2+\alpha^2\). Таким образом,
граничные условия для уравнения~\eqref{eq:3:3} при всех \(a\in U\) имеют
вид
\begin{equation}\label{eq:4:5}
	\begin{aligned}
		w(-\zeta_2)&=w(\zeta_2),\\
		w'(-\zeta_2)&=w'(\zeta_2)+\gamma w(\zeta_2),
	\end{aligned}
\end{equation}
и при этом
\begin{equation}\label{eq:4:6}
	\gamma=\left\{
	\begin{aligned}
		\zeta_2\left(\zeta_2^2-\alpha^2\right)^{-1}&&(a\in U_1\cup U_2)\\
		\zeta_2\left(\zeta_2^2+\alpha^2\right)^{-1}&&(a\in U_3).
	\end{aligned}\right.
\end{equation}
Чтобы получить уравнение для спектра, запишем общее решение
уравнения~\eqref{eq:3:3} в виде
\begin{equation}\label{eq:4:7}
	w(\zeta)=Aw_1(\zeta)+Bw_2(\zeta),
\end{equation}
где \(w_i(\zeta)\) (\(i=1,2\)) "--- решения, построенные в~\ref{pt:3}. Подставляя
выражение~\eqref{eq:4:7} в граничные условия~\eqref{eq:4:5}, получим систему
двух однородных уравнений для коэффициентов \(A\) и \(B\). Приравнивая определитель
этой системы нулю, имеем уравнение для спектра
\begin{equation}\label{eq:4:8}
	[w_1(\zeta_2)-w_1(-\zeta_2)]\,[w_2'(\zeta_2)-w_2'(-\zeta_2)+\gamma
	w_2(\zeta_2)]=[w_2(\zeta_2)-w_2(-\zeta_2)]\,
	[w_1'(\zeta_2)-w_1'(-\zeta_2)+\gamma w_1(\zeta_2)].
\end{equation}
Заметим, что полученная задача для уравнения~\eqref{eq:3:3} уже не является
периодической.

\section{Доказательство леммы о спектре}\label{pt:5}
Лемма о спектре, сформулированная в~\ref{pt:2}, содержит пять утверждений,
которые последовательно доказываются в этом параграфе.

\subsection{Доказательство утверждения~(1)}\label{pt:5:1}
В этом пункте предполагается, что \(a\in A_1\), т.~е. \(a_0\leqslant a<a_1\).

Введём сокращённые обозначения
\begin{align*}
	U&\equiv U\left(-\dfrac{1}{2}\lambda\alpha^2,
	\zeta_2\sqrt{2\lambda}\right),&
	U^{(1)}&\equiv U^{(1)}\left(-\dfrac{1}{2}\lambda\alpha^2,
	\zeta_2\sqrt{2\lambda}\right),\\
	\overline{U}&\equiv \overline{U}\left(-\dfrac{1}{2}\lambda\alpha^2,
	\zeta_2\sqrt{2\lambda}\right),&
	\overline{U}^{(1)}&\equiv \overline{U}^{(1)}\left(-\dfrac{1}{2}\lambda
	\alpha^2,\zeta_2\sqrt{2\lambda}\right).
\end{align*}
Из соотношений~\eqref{eq:3:13} следует, что
\begin{equation}\label{eq:5:1}
	\begin{aligned}
		w_2(\zeta_2)&=U+R_1(\zeta_2),& w_2(\zeta_2)&=\overline{U}+R_2(\zeta_2),\\
		w_2'(\zeta_2)&=\sqrt{2\lambda}\,U^{(1)}+R_2'(\zeta_2),&
		w_2'(\zeta_2)&=\sqrt{2\lambda}\,\overline{U}^{(1)}+R_2'(\zeta_2).
	\end{aligned}
\end{equation}
Для вычисления величин \(w_i(-\zeta_2)\), \(w_i'(-\zeta_2)\) (\(i=1,2\)) надо
воспользоваться определением~\eqref{eq:3:17} и уравнениями связи~\eqref{eq:3:18}.
Тогда
\begin{equation}\label{eq:5:2}
	\begin{aligned}
		w_1(-\zeta_2)&=A_1[U+R_3(-\zeta_2)]+B_1[\overline{U}+R_4(-\zeta_2)],\\
		w_2(-\zeta_2)&=A_2[U+R_3(-\zeta_2)]+B_2[\overline{U}+R_4(-\zeta_2)],\\
		w_1'(-\zeta_2)&=A_1[-\sqrt{2\lambda}U^{(1)}+R_3'(-\zeta_2)]+
		B_1[-\sqrt{2\lambda}\overline{U}^{(1)}+R_4'(-\zeta_2)],\\
		w_2'(-\zeta_2)&=A_2[-\sqrt{2\lambda}U^{(1)}+R_3'(-\zeta_2)]+
		B_2[-\sqrt{2\lambda}\overline{U}^{(1)}+R_4'(-\zeta_2)],
	\end{aligned}
\end{equation}
где коэффициенты \(A_i\), \(B_i\) удовлетворяют равенствам~\eqref{eq:3:19}.
Учитывая определение величин \(R_i(\zeta)\) (\(i=3,4\)), данное в~\ref{pt:3},
замечаем, что достаточно оценить величины \(R_i(\zeta_2)\), \(R'(\zeta_2)\)
(\(i=1,2\)), и полученные оценки будут справедливы для величин \(R_i(-\zeta_2)\),
\(R_i'(-\zeta_2)\) соответственно при \(i=3,4\). Общая для всех случаев
\(a\in A_i\) (\(i=1-5\)) схема доказательства соответствующих утверждений из леммы
о спектре состоит в подстановке выражений типа~\eqref{eq:5:1} и~\eqref{eq:5:2}
для \(w_i(\pm\zeta_2)\), \(w_i'(\pm\zeta_2)\) в уравнение~\eqref{eq:4:8}
для спектра и использования затем асимптотик функций Вебера. При этом
оказывается, что все нужные асимптотики известны и содержатся в работах
\cite{24,25}.

При \(a\in A_1\) начнём со случая, когда \(\lambda\alpha^2>
\lambda^{\varepsilon}\) (\(\varepsilon>0\)), и введём величины
\begin{equation}\label{eq:5:3}
	\mu^2=\lambda\alpha^2,\qquad t=\zeta_2\alpha^{-1},\qquad
	(\mu t\sqrt{2}=\zeta_2\sqrt{2\lambda}).
\end{equation}
Так как \(\varphi(\zeta_2)=\pi\), \(\varphi(x_2)=\alpha\) и \(\pi-x_2>C>0\),
то при \(a\in A_1\) выполняются неравенства \(t-1>C>0\). Тогда можно
воспользоваться асимптотическим разложением из работы \cite{24}, справедливым
при \(t>1\). Получаем разложение вида
\begin{equation}\label{eq:5:4}
	U=l_1(\mu)\dfrac{e^{-\mu^2\xi(t)}}{(t^2-a)^{1/4}}
	\sum\limits_{s=0}^{\infty}\dfrac{u_s(t)}{(t^2-1)^{3s/2}}\cdot
	\dfrac{1}{\mu^{2s}}.
\end{equation}
В этой формуле \(u_s(t)\) "--- полином степени \(3s\) при нечётных \(s\)
и степени \(3s-2\) "--- для чётных \(s\). В частности,
\begin{align*}
	u_0(t)&=1,& u_1(t)&=\dfrac{1}{24}(t^3-6t),&
	u_2(t)&=\dfrac{1}{1152}(9t^4+249t^2+145).
\end{align*}
Функция \(\xi(t)\) в~\eqref{eq:5:4} определяется равенством
\begin{equation}\label{eq:5:5}
	\xi(t)=\int\limits_1^t\sqrt{t^2-1}\,dt=\dfrac{t}{2}\sqrt{t^2-1}-
	\dfrac{1}{2}\ln\left(1+\sqrt{t^2-1}\right).
\end{equation}
Соответствующие разложения для производной \(U^{(1)}\) имеют вид (см. \cite{24})
\begin{equation}\label{eq:5:6}
	U^{(1)}=-\dfrac{\mu l_1(\mu)}{\sqrt{2}}(t^2-1)^{1/4}e^{-\mu^2\xi(t)}
	\sum\limits_{s=0}^{\infty}\dfrac{v_s(t)}{(t^2-1)^{3s/2}}\cdot
	\dfrac{1}{\mu^{2s}},
\end{equation}
\begin{align*}
	v_0(t)&=1,& v_1(t)&=\dfrac{1}{24}(t^3+6t),&
	v_2(t)&=\dfrac{1}{1152}(15t^4-327t^2-143),
\end{align*}
и степени полиномов \(v_s(t)\) такие же, как и полиномов \(u_s(t)\).
Аналогичные разложения имеют место и для величин \(\overline{U}\),
\(\overline{U}^{(1)}\):
\begin{equation}\label{eq:5:7}
	\begin{aligned}
		\overline{U}&=2l_1(\mu)\dfrac{e^{\mu^2\xi(t)}}{(t^2-1)^{1/4}}
		\sum\limits_{s=0}^{\infty}\dfrac{(-1)^su_s(t)}{(t^2-1)^{3s/2}}
		\cdot\dfrac{1}{\mu^{2s}},\\
		\overline{U}^{(1)}&=\sqrt{2}l_1(\mu)(t^2-1)^{1/4} e^{\mu^2\xi(t)}
		\sum\limits_{s=0}^{\infty}\dfrac{(-1)^sv_s(t)}{(t^2-1)^{3s/2}}
		\cdot\dfrac{1}{\mu^{2s}}.
	\end{aligned}
\end{equation}
Явный вид константы \(l_1(\mu)\) нам не важен, так как \(l_1(\mu)>C>0\)
при \(\mu^2\gg 1\) (см. \cite{24}). Как уже указывалось в~\ref{pt:2},
величина \(a_0\), входящая в определение области \(A_1\), выбрана таким
образом, что \(\pi-x_2>C\). Используя это и полагая \(\zeta=\zeta_2\)
в~\eqref{eq:3:5}, получим, что при \(a\in A_1\) выполняется \(t^2-1\geqslant
Ct^2\) и, следовательно, \(|u_s(t)(t^2-1)^{-3s/2}|\leqslant C_s\),
\(|v_s(t)(t^2-1)^{-3s/2}|\leqslant C_1\). Нужные нам асимптотики при
\(\mu^2>\lambda^{\varepsilon}\) (\(\varepsilon>0\)) теперь прямо следуют
из разложений~\eqref{eq:5:4}, \eqref{eq:5:6}, \eqref{eq:5:7} и имеют вид
\begin{equation}\label{eq:5:8}
	\begin{aligned}
		U&=l_1(\mu)e^{-\mu^2\xi(t)}(t^2-1)^{-1/4}\,
		\left[1+O(\mu^{-2})\right],\\
		U^{(1)}&=-\mu l_1(\mu)(t^2-1)^{1/4}e^{-\mu^2\xi(t)}\,
		\left[1+O(\mu^{-2})\right],\\
		\overline{U}&=2l_1(\mu)(t^2-1)^{-1/4}e^{\mu^2\xi(t)}\,
		\left[1+O(\mu^{-2})\right],\\
		U^{(1)}&=\sqrt{2}\mu l_1(\mu)(t^2-1)^{1/4}e^{\mu^2\xi(t)}\,
		\left[1+O(\mu^{-2})\right].
	\end{aligned}
\end{equation}
Так как в соответствии с определениями~\eqref{eq:5:3}, \eqref{eq:5:5}
и~\eqref{eq:2:2} выполняется \(\int_{x_2}^{\pi}\sqrt{f}\,dx=
\alpha^2\xi(t)\), то в области \(A_1\) при \(\mu^2>\lambda^{\varepsilon}\)
выполняются неравенства \(\xi(t)>C\alpha^{-2}\), \(\mu^2\xi(t)>C\lambda\),
\(t^2-1>C\), и тогда из~\eqref{eq:5:8} следует, что
\begin{equation}\label{eq:5:9}
	U=l_1(\mu)\,O\left(e^{-C\lambda}\right),\qquad
	U^{(1)}=l_1(\mu)\,O\left(e^{-C\lambda}\right).
\end{equation}
Оценим величины \(R_i(\zeta_2)\), \(R_i'(\zeta_2)\) (\(i=1,2\)). В интересующем
нас случае в равенствах~\eqref{eq:3:14} будет \(\beta=-\mu^2/2\),
\(z=\alpha t\sqrt{2\lambda}\), и из~\eqref{eq:3:15} следует, что
\(\rho(\beta)\simeq\alpha\sqrt{2\lambda}\). Таким образом, \(z>\rho(\beta)\),
и можно воспользоваться равенствами~\eqref{eq:3:16}, из которых с учётом
оценок~\eqref{eq:5:9} получим, что
\begin{equation}\label{eq:5:10}
	\begin{aligned}
		R_1(\zeta_2)&=l_1(\mu)\,O\left(e^{-C\lambda}\right),&
		R_1'(\zeta_2)&=l_1(\mu)\,O\left(e^{-C\lambda}\right),\\
		R_2(\zeta_2)&=\overline{U}\,O\left(\lambda^{-2/3}
		\ln\lambda\right),&
		R_2'(\zeta_2)&=\overline{U}^{(1)}\,O\left(\lambda^{-1/6}
		\ln\lambda\right).
	\end{aligned}
\end{equation}
Как уже указывалось, величины \(R_3(-\zeta_2)\), \(R_3'(-\zeta_2)\)
оцениваются так же как \(R_1(\zeta_2)\), \(R_1'(\zeta_2)\), а величины
\(R_4(-\zeta_2)\), \(R_4'(-\zeta_2)\) "--- как \(R_2(\zeta_2)\),
\(R_2'(\zeta_2)\). Используя эти оценки в равенствах~\eqref{eq:5:1}, имеем
\begin{equation}\label{eq:5:11}
	\begin{aligned}
		w_1(\zeta_2)&=l_1(\mu)\,O\left(e^{-C\lambda}\right),&
		w_1'(\zeta_2)&=l_1(\mu)\,O\left(e^{-C\lambda}\right),\\
		w_2(\zeta_2)&=\overline{U}\cdot\left(1+O\left(\lambda^{-2/3}
		\ln\lambda\right)\right),&
		w_2'(\zeta_2)&=\overline{U}^{(1)}\cdot\left(1\sqrt{2\lambda}+
		O\left(\lambda^{-1/6}\ln\lambda\right)\right).
	\end{aligned}
\end{equation}
Соответственно, из~\eqref{eq:5:2} получаем, что
\begin{equation}\label{eq:5:12}
	\begin{aligned}
		w_1(-\zeta_2)&=B_1\overline{U}\,\left(1+O\left(\lambda^{-2/3}
		\ln\lambda\right)\right)+l_1(\mu)\,O\left(e^{-C\lambda}\right),\\
		w_2(-\zeta_2)&=B_2\overline{U}\,\left(1+O\left(\lambda^{-2/3}
		\ln\lambda\right)\right)+l_1(\mu)\,O\left(e^{-C\lambda}\right),\\
		w_1'(-\zeta_2)&=B_1\overline{U}^{(1)}\,\left(-\sqrt{2\lambda}+
		O\left(\lambda^{-1/6}\ln\lambda\right)\right)+
		l_1(\mu)\,O\left(e^{-C\lambda}\right),\\
		w_2'(-\zeta_2)&=B_2\overline{U}^{(1)}\,\left(-\sqrt{2\lambda}+
		O\left(\lambda^{-1/6}\ln\lambda\right)\right)+
		l_1(\mu)\,O\left(e^{-C\lambda}\right).
	\end{aligned}
\end{equation}
Подставляя выражения~\eqref{eq:5:11}, \eqref{eq:5:12} в уравнение для
спектра~\eqref{eq:4:8} и используя асимптотики~\eqref{eq:5:8}, после простых
преобразований приводим уравнение для спектра при \(a\in A_1\) и \(\mu^2>
\lambda^{\varepsilon}\) к виду
\begin{equation}\label{eq:5:13}
	B_1\left[2+\dfrac{\gamma}{(t^2-1)^{1/4}\mu\sqrt{\lambda}}+
	\left(1+\dfrac{\gamma}{(t^2-1)^{1/4}\mu\sqrt{\lambda}}\right)\,
	O\left(\lambda^{-2/3}\ln\lambda+\mu^{-2}\right)\right]=
	O\left(e^{-C\lambda}\right).
\end{equation}
С учётом выражений~\eqref{eq:3:19} для \(B_1\) и~\eqref{eq:4:6} для \(\gamma\),
отсюда следует, что
\[
	\dfrac{\pi\lambda\alpha^2}{2}=n\pi+\dfrac{\pi}{2}+O\left(\lambda^{-2/3}
	\ln\lambda\right)\qquad(n\in\mathbb Z^+,\,n\gg 1).
\]
Остаётся воспользоваться выражением~\eqref{eq:2:1} для \(\alpha^2\) в области
\(U_1\), чтобы записать уравнение для спектра в виде
\[
	\lambda\int\limits_{x_1}^{x_2}\sqrt{g(x)-a}\,dx=n\pi+\dfrac{\pi}{2}+
	O\left(\lambda^{-2/3}\ln\lambda\right).
\]
Две ветви спектра в уравнении~\eqref{eq:2:12} получаются отсюда при \(n=2p\),
\(n=2p-1\), и утверждение~(1) при \(\mu^2\geqslant\lambda^{\varepsilon}\)
доказано.

Рассмотрим случай \(\mu^2\leqslant\lambda^{\varepsilon}\). В этом случае надо
воспользоваться асимптотическими разложениями
\begin{align*}
	U(\beta,z)&=e^{-z^2/4}z^{-\beta-1/2}\,\left[1-\dfrac{(\beta+1/2)
	(\beta+3/2)}{2z^2}+\dfrac{(\beta+1/2)(\beta+3/2)(\beta+5/2)
	(\beta+7/2)}{2\,(2z^2)^2}+\ldots\right],\\
	\overline{U}(\beta,z)&=\sqrt{\dfrac{2}{\pi}}\Gamma\left(\dfrac{1}{2}
	-\beta\right)e^{z^2/4}z^{\beta-1/2}\,\left[1+\dfrac{(\beta-1/2)
	(\beta-3/2)}{2z^2}+\dfrac{(\beta-1/2)(\beta-3/2)(\beta-5/2)
	(\beta-7/2)}{2\,(2z^2)^2}+\ldots\right]
\end{align*}
из работы \cite{25}. Эти разложения справедливы при \(z\gg|\beta|\), и их
можно дифференцировать по \(z\). При \(z\gg|\beta|\), таким образом, получаем
\begin{equation}\label{eq:5:14}
	\begin{aligned}
		U(\beta,z)&=e^{-z^2/4}z^{-\beta-1/2}\left[1+
		O\left(\dfrac{1+\beta^2}{z^2}\right)\right],\\
		U^{(1)}(\beta,z)&=-\dfrac{1}{2}e^{-z^2/4}z^{-\beta+1/2}\left[1+
		O\left(\dfrac{1+\beta^2}{z^2}\right)\right],\\
		\overline{U}(\beta,z)&=\sqrt{\dfrac{2}{\pi}}\Gamma\left(\dfrac{1}{2}-
		\beta\right)e^{z^2/4}z^{\beta-1/2}\left[1+
		O\left(\dfrac{1+\beta^2}{z^2}\right)\right],\\
		\overline{U}^{(1)}(\beta,z)&=\sqrt{\dfrac{1}{2\pi}}\Gamma\left(
		\dfrac{1}{2}-\beta\right)e^{z^2/4}z^{\beta+1/2}\left[1+
		O\left(\dfrac{1+\beta^2}{z^2}\right)\right].
	\end{aligned}
\end{equation}
Чтобы найти соответствующие асимптотики при \(z<0\), \(|z|\gg|\beta|\), надо
снова воспользоваться уравнениями связи~\eqref{eq:3:12}. В интересующем нас
случае \(-\beta=-\dfrac{1}{2}\lambda\alpha^2\), \(z=\zeta_2\sqrt{2\lambda}\).
Так как \(\zeta_2>C\), то \(|\beta z^{-1}|=O\left(\lambda^{1/2+\varepsilon}
\right)\), и при \(\varepsilon<1/2\) имеем \(|\beta z^{-1}|\ll 1\).
Определённая в~\ref{pt:3} функция \(\rho(\beta)\) \eqref{eq:3:15} возрастает
с ростом \(\beta\) и \(\beta(0)=0\) \cite{17}. Так как в рассматриваемом
случае \(|\beta|\leqslant C\lambda^{\varepsilon}\) и \(|z|\geqslant C
\lambda^{1/2}\), то, как и выше, \(z>\rho(\beta)\). Всё дальнейшее полностью
аналогично тому, что было сделано выше при \(\mu^2>\lambda^{\varepsilon}\),
только теперь используются асимптотики~\eqref{eq:5:14}. Уравнение для спектра
снова приводится к виду \(B_1=O\left(e^{-C\lambda}\right)\)
(см.~\eqref{eq:5:13}).

Утверждение~(1) доказано.

\subsection{Доказательство утверждения~(2)}\label{pt:5:2}
В этом пункте рассматривается область \(A_2\), в которой \(a_2+
C\lambda^{-1+\varepsilon_1}\leqslant a\leqslant a_0'\). Как и в
пункте~\ref{pt:5:1}, задача состоит в том, чтобы найти нужные асимптотики
величин \(w_i(\pm\zeta_2)\), \(w_i'(\pm\zeta_2)\) (\(i=1,2\)), где функции
\(w_i(\zeta)\) задаются равенствами~\eqref{eq:3:20} при \(\zeta>0\)
и равенствами~\eqref{eq:3:18}, \eqref{eq:3:23} при \(\zeta<0\). Величины
\(\mu\), \(t\), как и выше, определяются соотношениями~\eqref{eq:5:3}.

По определению~\eqref{eq:2:1} в области \(A_2\) выполняется равенство
\(\alpha^2=4\pi^{-1}\int_0^{x_2}\sqrt{a-g(x)}\,dx\). Так как \(x_{min}=0\),
то при малых \(x\) выполняется \(g(x)=a_2+\dfrac{1}{2}g^{(2)}(0)x^2+
O\left(x^4\right)\), \(g^{(2)}(0)>C>0\) и, следовательно, \(x_2\simeq
\left(2(a-a_2)/g_2(0)\right)^{1/2}\). Таким образом, \(\alpha^2\asymp
a-a_2\) и \(\mu^2>C\lambda^{\varepsilon_1}\) при \(a\in A_2\). Из условия
\(\pi-x_2>C\) следует, что \(\zeta_2\alpha^{-1}-1>C>0\) и \(t^2-1>Ct^2\).
При выполнении этих условий нужные асимптотики функций \(W\left(
\dfrac{\mu^2}{2},\pm\mu t\sqrt{2}\right)\) содержатся в работе \cite{24},
где показано, что при \(t>1\) имеют место разложения
\begin{equation}\label{eq:5:15}
	\begin{aligned}
		W\left(\dfrac{\mu^2}{2},\mu t\sqrt{2}\right)&=
		\dfrac{e^{-\pi\mu^2/4}l(\mu)}{\sqrt{2}(t^2-1)^{1/4}}\,
		\left(\cos\phi\,S_1-\sin\phi\,S_2\right),\\
		W\left(\dfrac{\mu^2}{2},\mu t\sqrt{2}\right)&=
		\dfrac{\sqrt{2}e^{\pi\mu^2/4}l(\mu)}{(t^2-1)^{1/4}}\,
		\left(\sin\phi\,S_1+\cos\phi\,S_2\right),
	\end{aligned}
\end{equation}
и в формулах~\eqref{eq:5:15}
\[
	S_1=\sum\limits_{s=0}^{\infty}\dfrac{(-1)^su_{2s}(t)}{(t^2-1)^{3s}
	\mu^{4s}},\qquad
	S_2=\sum\limits_{s=0}^{\infty}\dfrac{(-1)^su_{2s+1}(t)}{(t^2-1)^{3s+3/2}
	\mu^{4s+2}},
\]
а полиномы \(u_s(t)\) те же, что и в формуле~\eqref{eq:5:4}. Функция \(\phi\equiv
\phi(t,\mu)\) в правой части равенств~\eqref{eq:5:15} определяется равенством
\begin{equation}\label{eq:5:16}
	\phi=\mu^2\xi(t)+\dfrac{\pi}{4},
\end{equation}
в котором величина \(\xi(t)\) определена в~\eqref{eq:5:5}. Константа \(l(\mu)\)
при \(\mu\gg 1\) имеет асимптотику вида
\[
	l(\mu)=2^{1/4}\mu^{-1/4}\sum\limits_{s=0}^{\infty}l_s\mu^{-4s}\qquad
	\left(l_0=1,\,l_1=-\dfrac{1}{1152}\ldots\right).
\]
Из~\eqref{eq:5:15} следует, что
\begin{equation}\label{eq:5:17}
	\begin{aligned}
		W\left(\dfrac{\mu^2}{2},\mu t\sqrt{2}\right)&=
		\dfrac{e^{-\pi\mu^2/4}l(\mu)}{\sqrt{2}(t^2-1)^{1/4}}\,
		\left[\left(1-\dfrac{P_2}{\mu^4}\right)\cos\phi-
		\dfrac{P_1}{\mu^2}\sin\phi+O\left(\mu^{-6}\right)\right],\\
		W\left(\dfrac{\mu^2}{2},-\mu t\sqrt{2}\right)&=
		\dfrac{\sqrt{2}e^{\pi\mu^2/4}l(\mu)}{(t^2-1)^{1/4}}\,
		\left[\left(1-\dfrac{P_2}{\mu^4}\right)\sin\phi+
		\dfrac{P_1}{\mu^2}\cos\phi+O\left(\mu^{-6}\right)\right].
	\end{aligned}
\end{equation}
В этих формулах
\begin{equation}\label{eq:5:18}
	P_1=\dfrac{u_1(t)}{(t^2-1)^{3/2}}=\dfrac{1}{24}
	\dfrac{t^3-6t}{(t^2-1)^{3/2}},\qquad
	P_2=\dfrac{u_2(t)}{(t^2-1)^{3}}=
	\dfrac{-9t^4+249t^2+145}{1152(t^2-1)^3}.
\end{equation}
Разложения~\eqref{eq:5:15} можно дифференцировать, и таким образом получаются
нужные разложения производных
\begin{equation}\label{eq:5:19}
	\begin{aligned}
		W^{(1)}\left(\dfrac{\mu^2}{2},\mu t\sqrt{2}\right)&=
		-\dfrac{1}{2}e^{-\pi\mu^2/4}\mu l(\mu)(t^2-1)^{1/4}\,
		\left[\left(1-\dfrac{P_2^1}{\mu^4}\right)\sin\phi+
		\dfrac{P_1^1}{\mu^2}\cos\phi+O\left(\mu^{-6}\right)\right],\\
		W^{(1)}\left(\dfrac{\mu^2}{2},-\mu t\sqrt{2}\right)&=
		-e^{\pi\mu^2/4}\mu l(\mu)(t^2-1)^{1/4}\,
		\left[\left(1-\dfrac{P_2^1}{\mu^4}\right)\cos\phi-
		\dfrac{P_1^1}{\mu^2}\sin\phi+O\left(\mu^{-6}\right)\right].
	\end{aligned}
\end{equation}
В формулах~\eqref{eq:5:19}
\begin{equation}\label{eq:5:20}
	P_1^1=\dfrac{1}{24}\dfrac{t^3+6t}{(t^2+1)^{3/2}},\qquad
	P_2^1=\dfrac{-81t^4+249t^2-143}{1152(t^2-1)^3}.
\end{equation}
Оценим величины \(R_i(\zeta_2)\), \(R_i'(\zeta_2)\) (\(i=1,2\)). В рассматриваемой
области \(b=\dfrac{\mu^2}{2}>C\lambda^{\varepsilon_1}\), и, следовательно
(см.~\eqref{eq:2:9}),
\begin{align*}
	k(b)&=\dfrac{1}{2}e^{-\pi\mu^2/2}\left(1+O\left(e^{-\pi\mu^2}\right)\right),&
	k(b)^{-1}&=2e^{\pi\mu^2}\left(1+O\left(e^{-\pi\mu^2}\right)\right),\\
	k(b)^{1/2}&=\dfrac{1}{\sqrt{2}}e^{-\pi\mu^2/4}
	\left(1+O\left(e^{-\pi\mu^2}\right)\right),&
	k(b)^{-1/2}&=\sqrt{2}e^{\pi\mu^2/4}\left(1+O\left(e^{-\pi\mu^2}\right)\right).
\end{align*}
Так как в соответствии с~\eqref{eq:3:21} выполняется \(z=\zeta_2\sqrt{2}\lambda>
\sigma(b)\simeq\alpha\sqrt{2}\lambda\), то можно воспользоваться
равенством~\eqref{eq:3:14}, в котором \(\beta=b=\dfrac{\mu^2}{2}\), \(E(\beta,z)=
1\), и, используя выражения~\eqref{eq:5:17} и~\eqref{eq:5:19}, получим, что
\begin{equation}\label{eq:5:21}
	\begin{aligned}
		R_i(\zeta_2)&=l(\mu)\,O\left(\lambda^{-2/3}\ln\lambda\right)&
		&(i=1,2),\\
		R_i'(\zeta_2)&=\mu l(\mu)\,O\left(\lambda^{-1/6}\ln\lambda\right)&
		&(i=1,2).
	\end{aligned}
\end{equation}
Это позволяет получить нужные выражения для \(w_i(\zeta_2)\), \(w_i'(\zeta_2)\)
(\(i=1,2\)) по формулам~\eqref{eq:3:20} после подстановки в них
выражений~\eqref{eq:5:17} и~\eqref{eq:5:19}. Для получения соответствующих
формул для величин \(w_i(-\zeta_2)\), \(w_i'(-\zeta_2)\) (\(i=1,2\)) надо
воспользоваться соотношениями~\eqref{eq:3:18}, \eqref{eq:3:24}
и асимптотиками~\eqref{eq:5:17}, \eqref{eq:5:19}. При этом, в связи со
сказанным в~\ref{pt:3}, величины \(R_i(-\zeta_2)\), \(R_i'(-\zeta_2)\)
(\(i=3,4\)) оцениваются так же, как \(R_i(\zeta_2)\), \(R_i'(\zeta_2)\)
(\(i=1,2\)) (см.~\eqref{eq:5:21}).

Будем считать, что величины \(w_i(\pm\zeta_2)\), \(w_i'(\pm\zeta_2)\)
найдены. Переходим к решению уравнения для спектра~\eqref{eq:4:8}, которое
запишем в виде
\begin{equation}\label{eq:5:22}
	H_1\cdot H_3=H_2\cdot H_4.
\end{equation}
Здесь и ниже
\begin{equation}\label{eq:5:23}
	\begin{aligned}
		H_1&=w_1(-\zeta_2)-w_2(\zeta_2),&
		H_2&=w_2(-\zeta_2)-w_2(\zeta_2),\\
		H_3&=w_2'(-\zeta_2)-w_2'(\zeta_2)-\gamma w_2(\zeta_2),&
		H_4&=w_1'(-\zeta_2)-w_1'(\zeta_2)-\gamma w_1(\zeta_2).
	\end{aligned}
\end{equation}
Используя найденные выражения \(w_i(\pm\zeta_2)\), \(w_i'(\pm\zeta_2)\)
(\(i=1,2\)), получим, что
\begin{align*}
	H_1&=\dfrac{2e^{\pi\mu^2/2}l(\mu)}{(t^2-1)^{1/4}}\,
	\left[E_-+O\left(\lambda^{-2/3}\ln\lambda\right)\right],\\
	H_2&=-\dfrac{l(\mu)}{(t^2-1)^{1/4}}\,
	\left[E_-+O\left(\lambda^{-2/3}\ln\lambda\right)\right],\\
	H_3&=-\mu l(\mu)\sqrt{\lambda}(t^2-1)^{1/4}\,
	\left[E_++O\left(\lambda^{-2/3}\ln\lambda\right)\right]+
	\dfrac{\gamma l(\mu)}{(t^2-1)^{1/4}}\,
	\left[E_-+O\left(\lambda^{-2/3}\ln\lambda\right)\right],\\
	H_4&=-2\mu l(\mu)\sqrt{\lambda}e^{\pi\mu^2/2}(t^2-1)^{1/4}\,
	\left[E_++O\left(\lambda^{-2/3}\ln\lambda\right)\right],
\end{align*}
и в этих равенствах
\begin{align*}
	E_-&=\sin\phi+\dfrac{P_1}{\mu^2}\cos\phi-\dfrac{P_2}{\mu^4}\sin\phi+
	O\left(\mu^{-6}\right),\\
	E_+&=\cos\phi-\dfrac{P_1^1}{\mu^2}\sin\phi-\dfrac{P_2^1}{\mu^4}\cos\phi+
	O\left(\mu^{-6}\right).
\end{align*}
Это позволяет переписать уравнение~\eqref{eq:5:22} в виде
\begin{equation}\label{eq:5:24}
	2E_-E_+=-\dfrac{\gamma(E_-)^2}{\mu\sqrt{\lambda}(t^2-1)^{1/2}}+
	O\left(\lambda^{-2/3}\ln\lambda\right).
\end{equation}
Перемножим выражения в левой части уравнения~\eqref{eq:5:24}. При сделанных
предположениях первый член в правой части уравнения~\eqref{eq:5:24} оценивается
как \(O\left(\lambda^{-1}\right)\). Учитывая, что \(P_1+P_1^1=\dfrac{1}{12}
\zeta_2^3(\zeta_2^2-\alpha^2)^{-3/2}\), и используя определение~\eqref{eq:5:16}
величины \(\phi\), получим следующее уравнение для спектра
\begin{equation}\label{eq:5:25}
	2\mu^2\xi(t)=k\pi+\dfrac{\pi}{2}-\dfrac{1}{12}\dfrac{\zeta_2^3}{
	(\zeta_2^2-\alpha^2)^{3/2}}\cdot\dfrac{1}{\mu^2}+O\left(\mu^{-6}+
	\lambda^{-2/3}\ln\lambda\right).
\end{equation}
Остаётся воспользоваться равенством \(2\mu^2\xi(t)=\lambda F(a)\), чтобы
получить нужный результат~\eqref{eq:2:13}. Различный выбор знаков в нём
соответствует \(k=2p\) и \(k=2p-1\). Утверждение~(2) доказано.

\subsection{Доказательство утверждения~(3)}\label{pt:5:3}
Для рассматриваемой в этом пункте области \(A_3-\mu^2\leqslant
C\lambda^{\varepsilon_2}\) нас интересуют асимптотики величин \(W(\beta,
\pm z)\), \(W^{(1)}(\beta,\pm z)\) при \(\beta=b=\mu^2/2\), \(z=\zeta_2\sqrt{2}
\lambda\). Если \(\varepsilon_2<1/2\), то \(\beta z^{-1}=O\left(\lambda^{-1/2+
\varepsilon_2}\right)\), и нужные асимптотические разложения имеются в
работе \cite{25}. В рассматриваемой области при любом знаке \(\beta\) они
имеют вид
\begin{equation}\label{eq:5:26}
	\begin{aligned}
		W(\beta,z)&=\left(\dfrac{2k(\beta)}{z}\right)^{1/2}
		(T_1\cos\theta-T_2\sin\theta),\\
		W(\beta,-z)&=\left(\dfrac{2}{zk(\beta)}\right)^{1/2}
		(T_1\sin\theta+T_2\cos\theta),\\
	\end{aligned}
\end{equation}
и в этих формулах величины \(T_i\) задаются асимптотическими равенствами
\begin{align*}
	T_1&=p_0+\dfrac{q_2}{1! 2z^2}-\dfrac{p_4}{2!(2z^2)^2}-\dfrac{q_6}{3!
	(2z^2)^3}+\dfrac{p_8}{4!(2z^2)^4}+\dfrac{q_{10}}{5!(2z^2)^5}\ldots,\\
	T_2&=-\dfrac{p_2}{1! 2z^2}-\dfrac{q_4}{2!(2z^2)^2}+\dfrac{p_6}{3!
	(2z^2)^3}+\dfrac{q_8}{4!(2z^2)^4}-\dfrac{p_{10}}{5!(2z^2)^5}-
	\dfrac{q_{12}}{6!(2z^2)^6}\ldots,
\end{align*}
в которых коэффициенты \(p_n\equiv p_n(\beta)\), \(q_n\equiv q_n(\beta)\)
определяются из соотношений \(p_n+iq_n=\dfrac{\Gamma(n+1/2+i\beta)}{
\Gamma(n+i\beta)}\). Таким образом,
\[
	p_0=1,\qquad p_2=\dfrac{3}{4}-\beta^2,\qquad p_4=\dfrac{105}{16}-
	\dfrac{43}{2}\beta^2+\beta^4,\qquad q_2=2\beta,\qquad q_4=2\beta
	\left(\dfrac{11}{2}-4\beta^2\right),
\]
и имеют место оценки \(|p_n|\leqslant C_n(1+\beta^n)\), \(|q_n|\leqslant
C_n(1+\beta^n)\). Общий член в рядах \(T_i\) оценивается, как и \(C_n
(1+\beta^n)z^{-n}\), и в рассматриваемой области \(z\gg 1\), \(\beta
z^{-1}\ll 1\). Таким образом, из~\eqref{eq:5:26} получаем нужные
представления
\begin{equation}\label{eq:5:27}
	\begin{aligned}
		W(\beta,z)&=\left(\dfrac{2k(\beta)}{z}\right)^{1/2}
		\left[\left(1+\dfrac{\beta}{z^2}\right)\cos\theta-
		\dfrac{4\beta-3}{8z^4}\sin\theta+O\left(\dfrac{1+\beta^4}{z^4}
		\right)\right],\\
		W(\beta,-z)&=\left(\dfrac{2}{zk(\beta)}\right)^{1/2}
		\left[\left(1+\dfrac{\beta}{z^2}\right)\sin\theta-
		\dfrac{4\beta-3}{8z^4}\cos\theta+O\left(\dfrac{1+\beta^4}{z^4}
		\right)\right].
	\end{aligned}
\end{equation}
Во всех этих формулах
\begin{equation}\label{eq:5:28}
	\theta=\dfrac14z^2-\beta\ln z+\dfrac{\pi}{4}+\dfrac12 \arg\Gamma
	\left(\dfrac12+i\beta\right),\qquad\left(\beta=b=\dfrac{\mu^2}{2}\right).
\end{equation}
Соответствующие выражения для производных получаются дифференцированием
равенств~\eqref{eq:5:26} и имеют вид
\begin{equation}\label{eq:5:29}
	\begin{aligned}
		W^{(1)}(\beta,z)&=-\left(\dfrac{zk(\beta)}{2}\right)^{1/2}\,
		\left[\left(1-\dfrac{\beta}{z^2}\right)\,\sin\theta+
		\dfrac{4\beta^2+5}{8z^4}\cos\theta+O\left(\dfrac{1+\beta^4}{z^4}
		\right)\right],\\
		W^{(1)}(\beta,-z)&=-\left(\dfrac{2}{zk(\beta)}\right)^{1/2}\,
		\left[\left(1-\dfrac{\beta}{z^2}\right)\,\cos\theta-
		\dfrac{4\beta^2+5}{8z^4}\sin\theta+O\left(\dfrac{1+\beta^4}{z^4}
		\right)\right].
	\end{aligned}
\end{equation}
Остаётся оценить величины \(R_i(\zeta_2)\), \(R_i'(\zeta_2)\) (\(i=1,2\))
и \(R_i(-\zeta_2)\), \(R_i'(-\zeta_2)\) (\(i=3,4\)). При \(b\gg 1\) из
определения~\eqref{eq:2:9} функции \(k(b)\) следует, что
\[
	k(b)=\sqrt2-1+\pi b\left(1-\dfrac{1}{\sqrt2}\right)+O(b^2),\qquad
	k(b)\leqslant C.
\]
Так как \(z>\sigma(b)\), то, используя равенства~\eqref{eq:5:27},
\eqref{eq:5:29}, \eqref{eq:3:14} и~\eqref{eq:3:22}, получим, что
\begin{equation}\label{eq:5:30}
	R_i(\zeta_2)=O\left(\lambda^{-11/12}\ln\lambda\right),\qquad
	R_i(\zeta_2)=O\left(\lambda^{1/12}\ln\lambda\right),\qquad
	(i=1,2).
\end{equation}
Как указывалось выше, эти оценки верны и для величин \(R_i(-\zeta_2)\),
\(R'(-\zeta_2)\) при \(i=3,4\), соответственно. Равенства~\eqref{eq:5:27}
и~\eqref{eq:5:29} вместе с оценками~\eqref{eq:5:30} дают нужные выражения
для величин \(w_i(\zeta_2)\), \(w_i'(\zeta_2)\) (\(i=1,2\)) (см. \eqref{eq:3:20}).
Соответствующие выражения для \(w_i(-\zeta_2)\), \(w_i'(-\zeta_2)\) (\(i=1,2\))
получаются с помощью формул~\eqref{eq:3:18}, \eqref{eq:3:23}, \eqref{eq:3:24}.
Используя все эти выражения, получим, что в уравнении~\eqref{eq:5:22}
\begin{align*}
	H_1&=-\dfrac{1}{k(b)}\sqrt{\dfrac{2}{z}}\left(F_-+O_1\right),&
	H_2&=\sqrt{\dfrac{2}{z}}\left(F_-+O_2\right),\\
	H_3&=\sqrt{\lambda z}\left(F_++O_3\right),&
	H_4&=-k(b)^{-1}\left(F_++O_4\right),
\end{align*}
и в этих формулах
\begin{gather*}
	\begin{aligned}
		F_-&=A(k(b)\cos\theta-\sin\theta)+B(k(b)\sin\theta+\cos\theta),\\
		F_+&=A'(k(b)\sin\theta-\cos\theta)+B'(k(b)\cos\theta+\sin\theta),
	\end{aligned}\\
	A=1+\dfrac{\mu^2}{2z^2},\qquad A'=1-\dfrac{\mu^2}{2z^2},\qquad
	B=-\dfrac{\mu^4-3}{8z^2},\qquad B'=\dfrac{\mu^4+5}{8z^2},\qquad
	(z=\zeta_2\sqrt{2\lambda}),
\end{gather*}
а все величины \(O_i\) допускают оценку вида \(O_i=
O\left(\mu^8\lambda^{-2}+\lambda^{-2/3}\ln\lambda\right)\). Уравнение для
спектра~\eqref{eq:5:22} теперь записывается в виде
\[
	F_+F_-=O\left(\mu^8\lambda^{-8}+\lambda^{-2/3}\ln\lambda\right).
\]
Замечаем, что, если в выражениях \(F_{\pm}\) ограничиться главными членами,
то \(F_-\sim k(b)\cos\theta-\sin\theta\), \(F_+\sim k(b)\sin\theta-\cos\theta\).
Таким образом, величины \(F_-\), \(F_+\) не могут быть малы одновременно.
Рассуждая так же, как в конце пункта~\ref{pt:5:2}, получим, что, как и в
пункте~\ref{pt:5:2}, уравнение для спектра приближённо факторизуется и
записывается в виде двух уравнений
\begin{equation}\label{eq:5:31}
	F_{\pm}=O\left(\mu^8\lambda^{-8}+\lambda^{-2/3}\ln\lambda\right).
\end{equation}
Учитывая, что \(2\theta=\Psi(\lambda,a)+\pi/2\), где функция \(\Psi(\lambda,
a)\) определена в~\eqref{eq:2:11} из~\eqref{eq:5:31} и полученных выражений
для функций \(F_{\pm}\), уже легко получается нужный результат~\eqref{eq:2:14}
с соответствующим выбором знаков. Утверждение~(3) доказано.

\subsection{Доказательство утверждения~(4)}\label{pt:5:4}
В области \(A_4\subset U_3\) при вычислении величин \(w_i(\pm\zeta_2)\),
\(w_i'(\pm\zeta_2)\) (\(i=1,2\)) будем исходить из равенств~\eqref{eq:3:25}
и оценок~\eqref{eq:3:26}. Так как в этой области \(\mu^2=\lambda\alpha^2<
C\lambda^{\varepsilon_2'}\), то можно воспользоваться равенствами~\eqref{eq:5:27},
\eqref{eq:5:29}, в которых \(\beta=b=-\mu^2/2\). Из этих равенств с
помощью~\eqref{eq:3:26} получаем, что в рассматриваемой области \(R_i(\pm\zeta_2)=
O\left(\lambda^{-5/4}\ln\lambda\right)\), \(R_i'(\pm\zeta_2)=
O\left(\lambda^{-1/4}\ln\lambda\right)\) (\(i=1,2\)). Дальнейший ход доказательства
утверждения~(4) полностью аналогичен приведённому в пункте~\ref{pt:5:3}
доказательства утверждения~(3). Утверждение~(4) доказано.

\subsection{Доказательство утверждения~(5)}\label{pt:5:5}
Рассуждая так же, как в пункте~\ref{pt:5:2}, получим, что в области \(A_5\)
выполняются соотношения \(\mu^2=\lambda\alpha^2\geqslant C
\lambda^{\varepsilon_1'}\). Равенства~\eqref{eq:3:25}, в которых \(b=-\mu^2/2\),
показывают, что для вычисления величин \(w_i(\pm\zeta_2)\), \(w_i'(\pm\zeta_2)\)
(\(i=1,2\)) нужны асимптотики функций \(W(\beta,z)\) при \(\beta=-\mu^2/2<0\).
Так как в области \(U_3\) соотношения~\eqref{eq:3:25} справедливы при любом
знаке \(\zeta\), то, в отличие от предыдущего, где параметр \(t\) определялся
равенством~\eqref{eq:5:3}, теперь будем считать, что \(t\) при любом \(z\)
определяется равенством \(z=\mu t\sqrt{2}\).

Согласно работе \cite{24}, при \(\beta=-\mu^2/2<0\) равномерно по \(t\)
(\(|t|<\infty\)) имеют место следующие асимптотические равенства
\begin{equation}\label{eq:5:32}
	\begin{aligned}
		W\left(-\dfrac{\mu^2}{2},\mu t\sqrt{2}\right)&=
		\dfrac{l(\mu)}{(t^2+1)^{1/4}}\left[\overline{S}_1\cos
		\overline{\phi}-\overline{S}_2\sin\overline{\phi}\right],\\
		W^{(1)}\left(-\dfrac{\mu^2}{2},\mu t\sqrt{2}\right)&=
		-\dfrac{\mu}{\sqrt2}l(\mu)\,(t^2+1)\left[\overline{S}_3\sin
		\overline{\phi}+\overline{S}_4\cos\overline{\phi}\right].
	\end{aligned}
\end{equation}
В этих равенствах
\begin{gather*}
	\overline{\phi}\equiv\overline{\phi}(t)=\mu^2\overline{\xi}(t)+
	\dfrac{\pi}{4},\\
	\overline{\xi}(t)=\int\limits_0^t\sqrt{t^2+1}\,dt=\dfrac{t}{2}
	\sqrt{t^2+1}+\dfrac12\ln(t+\sqrt{t^2+1})
\end{gather*}
и \(\overline{S}_i\) "--- асимптотические ряды вида
\begin{align*}
	\overline{S}_1&=\sum\limits_{s=0}^{\infty}\dfrac{(-1)^s
	\overline{u}_{2s}(t)}{(t^2+1)^{3s}\mu^{4s}},&
	\overline{S}_2&=\sum\limits_{s=0}^{\infty}\dfrac{(-1)^s
	\overline{u}_{2s+1}(t)}{(t^2+1)^{3s+3/2}\mu^{4s+2}},\\
	\overline{S}_3&=\sum\limits_{s=0}^{\infty}\dfrac{(-1)^s
	\overline{v}_{2s}(t)}{(t^2+1)^{3s}\mu^{4s}},&
	\overline{S}_4&=\sum\limits_{s=0}^{\infty}\dfrac{(-1)^s
	\overline{v}_{2s+1}(t)}{(t^2+1)^{3s+3/2}\mu^{4s+2}}.
\end{align*}
Как и фигурирующие в равенствах~\eqref{eq:5:4}, \eqref{eq:5:6} величины,
\(u_s(t)\), \(v_s(t)\), \(\overline{u}_s(t)\) и \(\overline{v}_s(t)\) "---
полиномы степени \(3s\) для нечётных \(s\) и степени \(3s-2\)
(\(s\geqslant 2\)) для чётных \(s\). При этом
\begin{equation}\label{eq:5:33}
	\begin{gathered}
		\overline{u}_0=\overline{v}_0(t)=1,\qquad
		\overline{u}_1=-\dfrac{1}{24}(t^3+6t),\qquad
		\overline{v}_1=-\dfrac{1}{24}(t^3-6t),\\
		\overline{u}_2(t)=\dfrac{1}{1152}(9t^4+249t^2-145),\qquad
		\overline{v}_2(t)=\dfrac{1}{1152}(-15t^4-327t^2+143).
	\end{gathered}
\end{equation}
Из разложений~\eqref{eq:5:32} следует, что
\begin{equation}\label{eq:5:34}
	\begin{aligned}
		W\left(-\dfrac{\mu^2}{2},\mu t\sqrt2\right)&=
		\dfrac{l(\mu)}{(t^2+1)^{1/4}}\left[\left(1-
		\dfrac{\overline{u}_2(t)}{\mu^4(t^2+1)^3}\right)\cos\overline{\phi}-
		\dfrac{\overline{u}_1(t)}{\mu^2(t^2+1)^{3/2}}\sin\overline{\phi}+
		O\left(\mu^{-6}\right)\right],\\
		W^{(1)}\left(-\dfrac{\mu^2}{2},\mu t\sqrt2\right)&=
		-\dfrac{\mu l(\mu)}{\sqrt2}(t^2+1)^{1/4}\left[\left(1-
		\dfrac{\overline{v}_2(t)}{\mu^4(t^2+1)^3}\right)\sin\overline{\phi}+
		\dfrac{\overline{v}_1(t)}{\mu^2(t^2+1)^{3/2}}+
		O\left(\mu^{-6}\right)\right].
	\end{aligned}
\end{equation}
Используя вытекающие отсюда грубые оценки
\[
	W\left(-\dfrac{\mu^2}{2},\mu t\sqrt2\right)=
	\dfrac{l(\mu)}{(t^2+1)^{1/4}}\,O(1),\qquad
	W^{(1)}\left(-\dfrac{\mu^2}{2},\mu t\sqrt2\right)=
	\mu l(\mu)(t^2+1)^{1/4}\,O(1),
\]
с помощью формул~\eqref{eq:3:22}, \eqref{eq:3:23} получаем следующие оценки
остаточных членов
\begin{equation}\label{eq:5:35}
	R_i(\zeta)=\dfrac{l(\mu)\,
	O\left(\lambda^{-1}\ln\lambda\right)}{(t^2+1)^{1/4}},\qquad
	R_i'(\zeta)=\mu l(\mu)(t^2+1)^{1/4}\,
	O\left(\lambda^{-1}\ln\lambda\right),\qquad (i=1,2).
\end{equation}
Ниже будут использоваться следующие обозначения: \(\phi_i=\overline{\phi}(t_i)\)
(\(i=1,2\)), где величины \(t_i\) определяются из соотношений \(\sqrt{2}\mu t_i=
\zeta_i\sqrt{2\lambda}\). Так как \(\zeta_2=-\zeta_1\), то \(t_1=-t_2\), и везде
ниже \(t\equiv t_2\). Из равенств~\eqref{eq:3:25}, \eqref{eq:5:34}
и~\eqref{eq:5:35} получаем значения \(w_i(\pm\zeta_2)\), \(w_i'(\pm\zeta_2)\)
(\(i=1,2\)), а затем и значения величин \(H_i\), входящих в уравнение для
спектра~\eqref{eq:5:22}. Учитывая, что
\[
	\phi_1+\phi_2=\dfrac{\pi}{2},\qquad
	\dfrac{\gamma}{\mu\sqrt{x}(t^2+1)^{1/2}}=O\left(\lambda^{-1}\right),
\]
после соответствующих вычислений получаем, что
\begin{equation}\label{eq:5:36}
	\begin{aligned}
		H_1&=\dfrac{l(\mu)}{(t^2+1)^{1/4}}\left(\overline{E}_2-
		\overline{E}_1+O_1\right),&
		H_2&=\dfrac{l(\mu)}{(t^2+1)^{1/4}}\left(\overline{E}_1-
		\overline{E}_2+O_2\right),\\
		H_3&=\mu l(\mu)\sqrt{\lambda}(t^2+1)^{1/4}\left(\overline{E}_3-
		E_4+O_3\right),&
		H_4&=\mu l(\mu)\sqrt{\lambda}(t^2+1)^{1/4}\left(\overline{E}_3-
		\overline{E}_4+O_4\right),
	\end{aligned}
\end{equation}
и при этом справедливы оценки
\[
	O_i=O\left(\mu^{-6}+\lambda^{-1}\ln\lambda\right).
\]
Величины \(\overline{E}_i\) в равенствах имеют следующий вид:
\begin{align*}
	\overline{E}_1&=\left(1+\dfrac{\overline{P}_2}{\mu^4}\right)
	\sin\phi_2-\dfrac{\overline{P}_1}{\mu^2}\cos\phi_1,&
	\overline{P}_1&=-\dfrac{\overline{u}_1(t)}{(t^2+1)^{3/2}},\\
	\overline{E}_2&=\left(1+\dfrac{\overline{P}_2}{\mu^4}\right)
	\cos\phi_2+\dfrac{\overline{P}_1}{\mu^2}\sin\phi_1,&
	\overline{P}_2&=\dfrac{\overline{v}_1(t)}{(t^2+1)^{3/2}},\\
	\overline{E}_3&=\left(1+\dfrac{\overline{P}_2^1}{\mu^4}\right)
	\cos\phi_2-\dfrac{\overline{P}_1^1}{\mu^2}\sin\phi_1,&
	\overline{P}_1^1&=-\dfrac{\overline{u}_2(t)}{(t^2+1)^{3}},\\
	\overline{E}_4&=\left(1+\dfrac{\overline{P}_2^1}{\mu^4}\right)
	\sin\phi_2+\dfrac{\overline{P}_1^1}{\mu^2}\cos\phi_1,&
	\overline{P}_2^1&=-\dfrac{\overline{v}_2(t)}{(t^2+1)^{3}}.
\end{align*}
Подставляя выражения~\eqref{eq:5:36} для \(H_i\) в уравнение~\eqref{eq:5:22}
для спектра, получаем, что это уравнение можно записать в виде
\begin{equation}\label{eq:5:37}
	(\overline{E}_2-\overline{E}_1)(\overline{E}_3-\overline{E}_4)=
	O\left(\mu^{-6}+\lambda^{-1}\ln\lambda\right).
\end{equation}
В отличие от предыдущих случаев, главные члены асимптотических разложений
величин \(\overline{E}_2-\overline{E}_1\) и \(\overline{E}_3-\overline{E}_4\)
совпадают и равны \(\cos\phi_2-\sin\phi_2\). Положим \(\phi_2=\pi/4+\pi p+
\delta\). Так как \(\overline{P}_1^1=-P_1+O\left(\lambda^{-1}\right)\), то
из~\eqref{eq:5:37} получаем следующее уравнение для \(\delta\):
\begin{equation}\label{eq:5:38}
	\sin^2\delta\left(1+\dfrac{\overline{P}_2+\overline{P}_2^1}{\mu^4}\right)-
	2\sin\delta\cos\delta\dfrac{\overline{P}_1}{\mu^2}+\cos^2\delta
	\dfrac{\left(\overline{P}_1\right)}{\mu^4}=
	O\left(\mu^{-6}+\lambda^{-1}\ln\lambda\right).
\end{equation}
Так как \(\overline{P}_1=(24)^{-1}+O\left(\lambda^{-2}\right)\),
\(|\overline{P}_2^1+\overline{P}_2^1\leqslant C\alpha^2\), то из~\eqref{eq:5:38}
получаем, что это уравнение имеет двойной корень вида \(\delta=-(48b)^{-1}+
O\left(|b|^{-3/2}+\lambda^{-1/2}(\ln\lambda)^{1/2}\right)\). Учитывая равенства
\[
	2\phi_2=2\mu^2\overline{\xi}(t)+\dfrac{\pi}{2},\qquad
	2\mu^2\overline{\xi}(t)=2\lambda\int\limits_0^{\pi}
	\sqrt{g(x)-a}\,dx,
\]
получаем отсюда нужный результат~\eqref{eq:2:16}.

Утверждение~(5) доказано, что и завершает доказательство леммы о спектре.

\section{Доказательство теоремы}\label{pt:6}
Результат леммы о спектре означает, что все точки положительного спектра,
если они существуют, асимптотически близки к решениям
уравнений~\eqref{eq:2:12}--\eqref{eq:2:16}.

Разобьём доказательство на три пункта~\ref{pt:6:1}, \ref{pt:6:2}, \ref{pt:6:3}
и докажем их последовательно.

\subsection{Доказательство существования двух ветвей спектра}\label{pt:6:1}
В этом пункте будет доказано, что во всём интервале \(U\) уравнения для
спектра~\eqref{eq:2:12}--\eqref{eq:2:16} можно записать в виде
\begin{equation}\label{eq:6:1}
	\lambda F(a)=2\pi p+H_{\pm}(b_2)+R_{\pm}(\lambda,a).
\end{equation}
При этом выполняются оценки
\begin{equation}\label{eq:6:2}
	\begin{aligned}
		|R_{\pm}(\lambda,a)&=O\left(\lambda^{-2/3}\ln\lambda\right),
		\qquad (a\geqslant a_2),\\
		|R_{\pm}(\lambda,a)&=O\left(\lambda^{-1/2}(\ln\lambda)^{1/2}\right),
		\qquad (a\leqslant a_2)
	\end{aligned}
\end{equation}
и решения уравнений~\eqref{eq:6:1} существуют.

Определим границу между областями \(A_2\) и \(A_3\) из условия \(b_2^3=
\lambda^{2/3}\). Так как при \(a\in A_1\cup A_2\) выполняется равенство
\(b_2=1/2\cdot\lambda\alpha_2^2\) и величина \(\alpha_2\equiv\alpha_2(a)\)
растёт с ростом \(a\), то в рассматриваемой области \(b_2\geqslant\lambda^{2/9}
\gg 1\). Тогда из~\eqref{eq:2:17} следует, что \(H_{\pm}(b_2)=\pm\pi/2-
(24b_2)^{-1}+O\left(b_2^{-3}\right)\), и нужный результат~\eqref{eq:6:1}
при \(a\in A_1\cup A_2\) получается из равенств~\eqref{eq:2:12} и~\eqref{eq:2:13}.

Рассмотрим области \(A_3\subset U_1\), в которой \(b_2\leqslant\lambda^{2/3}\).
В этой области, в силу первого из равенств~\eqref{eq:2:10},
\[
	F(a)=2\int\limits_{x_2}^{\pi}\sqrt{g(x)-a}\,dx=\zeta_2\sqrt{\zeta_2^2-
	\alpha^2}-\alpha^2\ln\left(\dfrac{\zeta_2}{\alpha}+\dfrac{1}{\alpha}
	\sqrt{\zeta_2^2-\alpha^2}\right).
\]
Учитывая определение~\eqref{eq:2:11} функции \(\Psi\) и определение~\eqref{eq:2:7}
функций \(H_{\pm}\), из уравнений~\eqref{eq:2:14} получим, что в рассматриваемой
области
\[
	\lambda F(a)=2\pi p+H_{\pm}(b_2)+O\left(\dfrac{b^4}{\lambda^2}+
	\lambda^{-2/3}\ln\lambda\right)
\]
и, так как при \(a\in A_3\) выполняется неравенство \(b_2\leqslant\lambda^{2/9}\),
то для \(a\geqslant a_2\) имеют место уравнения~\eqref{eq:6:1} и оценки~\eqref{eq:6:2}.

При \(a\leqslant a_2\) доказательство аналогично приведённому, надо только выбрать
границу между областями \(A_4\) и \(A_5\) из условия \(|b|=\lambda^{1/3}\).

Остаётся доказать, чо решения \(\lambda_{\pm}(a,p)\) уравнений~\eqref{eq:6:1}
существуют. Из результатов Олвера следует, что все остаточные члены,
обозначаемые в~\ref{pt:3} через \(R_i(\lambda)\), являются непрерывными
функциями \(\lambda\) (\(\lambda\gg 1\)). Отсюда следует, что и функции
\(R_{\pm}(\lambda,a)\) в уравнениях~\eqref{eq:6:1} являются непрерывными
функциями \(\lambda\) при \(\lambda\gg 1\).

Рассмотрим случай \(a\in U_1\). Так как в этой области \(b_2\geqslant C\lambda\)
(\(\alpha_2\geqslant C\)), то уравнение~\eqref{eq:6:1} совпадает с
уравнением~\eqref{eq:2:12}. Запишем это уравнение в виде \(f(x)=0\), где
\[
	f(x)=x-x_0^{\pm}-R_{\pm}\left(F^{-1}x,a\right),\qquad
	x=\lambda F,\qquad x_0^{\pm}=2\pi i\pm\dfrac{\pi}{2}.
\]
Так как функция \(f(x)\) непрерывна и меняет знак внутри отрезка
\([x_0^{\pm}-\Delta x,x_0^{\pm}+\Delta x]\), \(\Delta x=O\left(F^{-1}
x_0^{\pm}\right)\), то на этом отрезке имеется корень \(x^{\pm}\) уравнения
\(f(x)=0\), причём \(x^{\pm}=x_0^{\pm}+R_{\pm}\left(F^{-1}x_0^{\pm}\right)+
R_{\pm}\left(x^{\pm}F^{-1}\right)-R_{\pm}\left(F^{-1}x_0^{\pm}\right)\).
Из оценки \(\Delta x\) следует, что \(\left|R_{\pm}\left(xF^{-1}\right)-
R_{\pm}\left(x_0^{\pm}F^{-1}\right)\right|=O\left((x_0^{\pm}/F)^{-2/3}
\ln(x^{\pm}/F)\right)\), и утверждение~(1) теоремы при \(a\in U_1\)
доказано.

Рассмотрим случай \(a\in U_2\cup U_3\). В этом случае \(F(a)\geqslant
C>0\), и уравнение~\eqref{eq:6:1} имеет вид \(f(\lambda)=0\), где
\[
	f(\lambda)=\lambda-\lambda_p^0-\dfrac{1}{F(a)}H_{\pm}(b_2)-
	\dfrac{1}{F_a}R_{\pm}(\lambda).
\]
Так как \(|H_{\pm}(b_2)|\leqslant C\) и \(|F|\geqslant C\), то на отрезке
\(|\lambda-\lambda_p^0|\leqslant C\) имеется корень \(\lambda=\lambda_p^0+
\Delta\lambda\) уравнения \(f(\lambda)=0\) и \(|\Delta\lambda|\leqslant C\).
При этом
\begin{equation}\label{eq:6:3}
	\begin{gathered}
		\Delta\lambda=\bigl(F(a)\bigr)^{-1}H_{\pm}\left(
		\dfrac{\lambda_p^0\alpha^2}{2}\right)-\bigl(F(a)\bigr)^{-1}
		\Delta H_{\pm}(\lambda)+\bigl(F(a)\bigr)^{-1}R_{\pm}(\lambda),\\
		\Delta H_{\pm}=H_{\pm}\left(\left(\lambda_p^0+\Delta\lambda\right)
		\dfrac{\alpha^2}{2}\right)-H_{\pm}\left(
		\dfrac{\lambda_p^0\alpha^2}{2}\right).
	\end{gathered}
\end{equation}
Так как \(|\Delta\lambda|\leqslant C\), то в оценке \(R_{\pm}(\lambda)\)
можно заменить \(\lambda\) на \(\lambda_p^0\). Покажем, что
\begin{equation}\label{eq:6:4}
	|\Delta H_{\pm}(\lambda)|\leqslant C\lambda^{-2/3}.
\end{equation}
Рассмотрим область \(\alpha^2\geqslant C\left(\lambda_p^0\right)^{-1+2/15}\),
в которой \(|b|\geqslant C\left(\lambda_p^0\right)^{2/15}\) (\(|b|^5\geqslant
C\left(\lambda_p^0\right)^{2/3}\)). Тогда оценка~\eqref{eq:6:4} следует из
асимптотических разложений~\eqref{eq:2:17}.

Рассмотрим область \(\alpha^2\leqslant C\left(\lambda_p^0\right)^{-1+2/15}\)
и будем исходить из оценок
\begin{equation}\label{eq:6:5}
	|H_{\pm}(x)|\leqslant C\max\left|\dfrac{dH_{\pm}}{dx}\right|\,
	|\Delta\lambda|\,\alpha^2,\qquad |\Delta\lambda|\leqslant C.
\end{equation}
Из определения~\eqref{eq:2:7} функций \(H_{\pm}(x)\) следует, что
\[
	\dfrac{dH_{\pm}}{dx}=\pm\dfrac{\pi e^{\pi x}}{1+e^{2\pi x}}+
	e_n|x|-\Im\left(i\psi\left(\dfrac{1}{2}+ix\right)\right)\qquad
	\left(\psi(x)=\dfrac{\Gamma'(x)}{\Gamma(x)}\right).
\]
Для оценки последнего члена в правой части этого равенства вопсользуемся
интегральным представлением для функции \(\psi(x)\) (см.~\cite{22})
\begin{equation}\label{eq:6:6}
	\psi(\tau)=\ln\tau-\dfrac{1}{2\tau}-2\int\limits_0^{\infty}
	t(t^2+\tau^2)^{-1}\,\left(e^{2\pi t}-1\right)^{-1}\,dt.
\end{equation}
Нас интересует случай \(\tau=1/2+ix\) и, следовательно, \(|\tau|\geqslant
1/2\). При \(|\tau|\geqslant 1/2\) из~\eqref{eq:6:6} получаем, что
\(|\psi(\tau)|\leqslant|\ln(|\tau|)|+C\), и, таким образом, \(\left|
\dfrac{dH_{\pm}}{dx}\right|\leqslant C+2|\ln(|x|)|\), а потому
оценка~\eqref{eq:6:4} следует из оценок~\eqref{eq:6:5}. Существование решений
\(\lambda_{\pm}(a,p)\) с оценками~\eqref{eq:2:6} доказано.

\subsection{Структура спектра}\label{pt:6:2}
В этом пункте будет доказана асимптотическая полнота двух построенных ветвей
спектра \(\lambda_{\pm}(a,p)\), т.~е. доказано, что весь положительный спектр
задачи~\eqref{eq:1:1} при \(\lambda\gg 1\) исчерпывается ветвями
\(\lambda_{\pm}(a,p)\).

При \(a_2-a>C\) задача дефинитна и структура спектра известна \cite{3,4,5}.
В этом случае ветви \(\lambda_{\pm}(a,p)\) действительно исчерпывают весь
спектр, так как он состоит из двух ветвей с одинаковой асимптотикой
(см. \cite{1,4,5}). При \(a>a_2\) задача индефинитна и основная трудность
состоит в том, что для периодической индефинитной задачи осцилляционная
теорема, по-видимому, не доказана. Условие чётности функции \(g(x)\)
позволяет обойти эту трудность и свести периодическую задачу к двум задачам
с разделёнными граничными условиями, а затем воспользоваться результатами
работы~\cite{9} (см.~также \cite{8,14}).

Рассмотрим две вспомогательные задачи для уравнения~\eqref{eq:1:1}
на интервале \((0,\pi)\) "--- задачу \(D\) (задачу Дирихле)
\begin{equation}\label{eq:6:7}
	y''+\lambda^2(g(x)-a)y=0,\qquad y(0)=y(\pi)=0,
\end{equation}
и задачу \(N\) (задачу Неймана)
\begin{equation}\label{eq:6:8}
	y''+\lambda^2(g(x)-a)y=0,\qquad y'(0)=y'(\pi)=0.
\end{equation}
Обозначим спектры этих задач через \(\{\lambda^D\}\) и \(\{\lambda^N\}\)
соответственно. Из условия чётности функции \(g\) следует, что спектр
\(\{\lambda\}\) задачи~\eqref{eq:1:1} является объединением спектров
\(\{\lambda^D\}\) и \(\{\lambda^N\}\): \(\{\lambda\}=\{\lambda^D\}\cup
\{\lambda^N\}\). При \(\lambda=\lambda^D\) собственная функция
задачи~\eqref{eq:1:1} является нечётным продолжением на интервал
\((-\pi,\pi)\) соответствующей собственной функции задачи \(D\)
\eqref{eq:6:7}, а при \(\lambda=\lambda^N\) "--- чётным продолжением
соответствующей собственной функции задачи \(N\) \eqref{eq:6:8}.

Заметим, что, хотя задачи \(D\) и \(N\) "--- задачи с одной точкой
поворота, нас интересует случай, когда точки поворота близки к концам
интервала \((0,\pi)\), в которых \(g'(0)=g'(\pi)=0\). В силу этого
при анализе задач \(D\) и \(N\) мы вновь используем замену~\eqref{eq:3:2}.

Асимптотики собственных значений \(\lambda^D\) и \(\lambda^N\)
строятся точно так же, как было сделано выше для задачи~\eqref{eq:1:1}.
Проведя это построение, получаем, что имеют место аналоги
формул~\eqref{eq:2:4}:
\begin{equation}\label{eq:6:9}
	\begin{aligned}
		\lambda^D(a,p)&=\lambda_p^0+F(a)^{-1}H_-(b_2(\lambda_p^0))+
		R_D(a,p),\\
		\lambda^N(a,p)&=\lambda_p^0+F(a)^{-1}H_+(b_2(\lambda_p^0))+
		R_N(a,p),
	\end{aligned}
\end{equation}
и, следовательно, \(\lambda^D(a,p)=\lambda_-(a,p)\), \(\lambda^N(a,p)=
\lambda_+(a,p)\). При этом для величин \(R_D\) и \(R_N\) справедливы
оценки~\eqref{eq:2:6}.

Докажем асимптотическую полноту спектров \(\{\lambda^D(a,p)\}\) и
\(\{\lambda^N(a,p)\}\), т.~е. отсутствие при \(p\gg 1\) других
собственных значений, близких к \(\{\lambda^D(a,p)\}\) и
\(\{\lambda^N(a,p)\}\). Это следует из осцилляционной теоремы, доказанной
в индефинитном случае для задачи Штурма--Лиувилля с разделёнными граничными
условиями в работе \cite{9}. Из результатов этой работы следует, что для
доказательства асимптотической полноты спектров \(\{\lambda^D(a,p)\}\)
и \(\{\lambda^N(a,p)\}\) достаточно доказать, что \(\forall k\gg 1\)
в задаче \(D\) (\(N\)) существует собственная функция, отвечающая
собственному значению \(\lambda^D(a,p)\) (\(\lambda^N(a,p)\)) для некоторого
\(p\gg 1\), и имеющая \(k\) нулей в интервале \((0,\pi)\). С помощью оценок
из работы \cite{17} задача сводится к вопросу о количестве нулей в интервале
\((0,\zeta_2)\) соответствующей функции Вебера. Используя результаты о
нулях этих функций, содержащиеся в работах \cite{17,24,25}, получим, что
функция \(w(\zeta)\) (\(u(x)=\left(\varphi'(\zeta)\right)^{1/2}w(\zeta)\)),
отвечающая собственному значению \(\lambda_-(a,p)\) в задаче \(D\), имеет
\(p-1\) нулей в интервале \((0,\zeta_2)\), а функция \(w(\zeta)\), отвечающая
собственному значению \(\lambda_+(a,p)\) в задаче \(N\), имеет \(p\) нулей
в интервале \((0,\zeta_2)\).

Таким образом, спектры \(\{\lambda^D(a,p)\}\) задачи \(D\) и
\(\{\lambda^N(a,p)\}\) задачи \(N\) асимптотически полны, а следовательно,
асимптотически полон и спектр \(\{\lambda_{\pm}(a,p)\}\) задачи~\eqref{eq:1:1}.
Утверждение~(1) теоремы доказано. Кроме того, доказано, что собственная
функция периодической задачи, отвечающая собственному значению
\(\lambda_{\pm}(a,p)\), имеет \(2p\) нулей на полуинтервале \([-\pi,\pi)\).

Заметим, что \(\lambda=0\) "--- собственное значение задачи \(N\), и
соответствующая собственная функция равна константе. Таким образом, естественно
предположить, что при \(a>a_2\), \(n\gg 1\), выполняются равенства~\eqref{eq:2:8}.
Это будет доказано в пункте~\ref{pt:6:3} независимо от результатов работы
\cite{9}.

\subsection{Доказательство второго утверждения теоремы}\label{pt:6:3}
Рассмотрим область \(U'_3\subset U_3\), в которой \(a_2-a>C>0\). В этой области
задача дефинитна. Переходя с помощью преобразования Лиувилля, описанного
в~\ref{pt:1}, к классической задаче Штурма--Лиувилля~\eqref{eq:1:2}
(см.~\cite{4}), снова получим, что собственная функция, отвечающая собственному
значению \(\lambda_{\pm}(a,p)\), имеет \(2p\) нулей на полуинтервале \([-\pi,
\pi)\). Тогда по осцилляционной теореме для периодической задаче Штурма--Лиувилля
\cite{3,4} получаем, что при \(a\in U_3'\), \(n\gg 1\), имеют место
равенства~\eqref{eq:2:8}. Покажем, что эти равенства справедливы \(\forall a\in
U\). Рассмотрим величины \(s_{n,a}\), \(s_{\pm}(a,p)\), и отрезки \(I_p\) вида
\begin{gather*}
	s_{n,a}=F(a)\lambda_{n,a},\qquad s_{\pm}(a,p)=F(a)\lambda_{\pm}(a,p),\\
	I_p=\left[2\pi p-\dfrac{\pi}{2}-\delta,2\pi p+\dfrac{\pi}{2}+\delta\right],
	\qquad\delta=\dfrac{\pi}{10}.
\end{gather*}
Интервалы \(I_p\) не перекрываются, и, так как \(|H_{\pm}(x)|\leqslant
\pi/2+\pi/20\) (см.~\ref{pt:2}), то при \(p\gg 1\) имеются точно два
собственных значения \(s_{\pm}(a,p)\in I_p\), \(\forall a\in U\).
Так как для любого фиксированного \(a\) при \(k\gg 1\), \(p\gg 1\)
имеет место взаимно однозначное соответствие \(\{s_{n,a}\}\leftrightarrow
\{s_{\pm}(a,p)\}\) и \(s_{2k-1,a},s_{2k,a}\in I_k\) при \(a\in U_{\rho}'\),
то утверждения~\eqref{eq:2:8} справедливы \(\forall a\in U\). При \(a>a_2\)
выполняется неравенство \(s_+(a,k)>s_-(a,k)\), и по определению порядка
\(s_{2k,a}\geqslant s_{2k-1,a}\). Отсюда следует, что \(s_{2k-1,a}=
s_-(a,k)\), \(s_{2k,a}=s_+(a,k)\) \(\forall k\gg 1\) и \(a>a_2\).
Утверждение~(2) доказано, что и завершает доказательство теоремы.

Полученные результаты показывают, что асимптотическая структура спектра
задачи~\eqref{eq:1:1} при \(g\in\mathcal G\), \(a_1-a>C\) не меняется
при переходе от дефинитной при \(a<a_2\) к индефинитной при \(a>a_2\)
задаче. При этом происходит только снятие асимптотического вырождения,
имеющего место при \(a<a_2\) (см.~\eqref{eq:1:5}). Остаётся справедливой
и осцилляционная теорема.

В заключение сделаем два замечания. Первое состоит в том, что полученные
результаты непосредственно обобщаются на периодическую задачу Штурма--Лиувилля
для уравнения вида \(y''+\lambda^2f(x,a)y=0\), если для этого уравнения
справедливы результаты Олвера. Отметим также, что приведённые в пункте~(1)
теоремы результаты при \(a-a_2>C\) справедливы и без предположения чётности.
Второе замечание состоит в том, что функции \(H_{\pm}(x)\), по-видимому,
описывают асимптотику положительного дискретного спектра в окрестности
точек типа \(a=a_2\), где рождаются (умирают) две простые точки поворота,
для зависящих от параметра \(a\) задач Штурма--Лиувилля, вне зависимости
от выбора граничных условий (см.~\eqref{eq:6:9}).

\end{document}